\documentclass{article}
\usepackage{latexsym}
\usepackage{mathrsfs}
\usepackage{indentfirst}
\usepackage{amsmath,amsthm, amssymb}
\usepackage{graphicx}
\usepackage{setspace}
\usepackage{setspace}
\usepackage{geometry}
\usepackage{float}
\usepackage{epsfig}
\usepackage{booktabs}
\usepackage{verbatim}
\usepackage{multirow}
\usepackage{tikz}
\usepackage{tikz-cd}
\usepackage{csquotes}
\usepackage{braids}
\usepackage[hidelinks]{hyperref}
\usepackage{bbm}
\newtheorem{theorem}{Theorem}[section]
\newtheorem{proposition}[theorem]{Proposition}
\newtheorem{lemma}[theorem]{Lemma}
\newtheorem{corollary}[theorem]{Corollary}

\newtheorem{definition}[theorem]{Definition}
\usepackage{caption}
\begin{document}
\title{Rationally Null-homologous Floer Simple Knot}
\author{Zhaojun Chen$^1$\thanks{Email address: zchen5@caltech.edu.}
\\ $^1$Department of Mathematics,
\\ California Institute of Technology, Pasadena, CA 91125, United States
}
\date{Sep 5, 2026}
\maketitle
%removing the first page number
%\thispagestyle{empty}

%abstract part
~\\
\begin{abstract}
We prove that every rationally null-homologous Floer simple knot in a closed, oriented, irreducible three-manifold with nonzero Thurston norm is unknotted. The proof combines a cut-and-reglue comparison for extremal Floer groups with Binns’s theorem on Floer-minimal knots in balanced sutured manifolds.
%key words
$\newline$\noindent\textbf{Key words}: Floer Simple Knot, Sutured Floer Homology, Broken Fibration, Surgery Exact Triangle.
\end{abstract}
~\\
\section{Introduction}
For a rationally null-homologous knot $K$ in a closed oriented
three-manifold $Y$, the knot Floer spectral sequence in [10] gives
\[
\operatorname{rank}\widehat{HFK}(Y,K)
\geq
\operatorname{rank}\widehat{HF}(Y).
\]
The knot is called \emph{Floer simple} when equality holds.

In [8], Ni proved that a null-homotopic Floer simple knot is unknotted whenever its ambient closed three-manifold is irreducible and has nonzero Thurston norm. He also proved the corresponding result for rationally null-homologous knots in surface bundles over \(S^1\) with fiber of genus greater than one. These results led him to ask whether a closed irreducible three-manifold with nonzero Thurston norm can contain a nontrivial rationally null-homologous Floer simple knot. 

Ni's results for null-homotopic Floer simple knot suggests a negative answer to the question above, but a direct extension of Ni’s argument encounters a difficulty in the cutting-and-regluing process. His proof uses an internal hierarchy to reduce to the surface-bundle case. The null-homotopy hypothesis provides control of the knot after cutting along incompressible surfaces and is used to arrange decomposing surfaces disjoint from the knot. Under the weaker assumption of being rationally null-homologous in the original manifold, the reglued knot need not remain rationally null-homologous. Consequently, the usual argument for transferring Floer simplicity cannot be applied unchanged to the intermediate pairs.

Binns’s work [2] on Floer-minimal knots in balanced sutured manifolds suggests a different route to the topological conclusion. Let \(K\) be a knot in the interior of a balanced sutured manifold \((N,\gamma)\), and let \((N(K),\gamma(K))\) denote $N\setminus\!\setminus K$ equipped with two additional meridional sutures. Suppose that $[K]=0$ in $H_1(N;\mathbb Q)/H_1(\partial N;\mathbb Q). $ Binns proves that the equality
$$ \operatorname{rank}SFH(N(K),\gamma(K)) = 2\operatorname{rank}SFH(N,\gamma) \neq 0 $$
forces \(N\) to have a closed connected summand containing \(K\), in which \(K\) is Floer simple. In particular, when \(N\) is irreducible and has nonempty boundary, this summand must be \(S^3\); unknot detection in \(S^3\) then implies that \(K\) is unknotted.

Crucially, the homological hypothesis in Binns’s theorem requires only that the rational homology class of \(K\) come from \(\partial N\), rather than that \(K\) be null-homotopic in $N$. This suggests replacing the reduction to a surface bundle by a construction of suitable balanced sutured manifolds to which Binns’s theorem applies. The central task is then to obtain the sutured rank equality from the Floer-simple hypothesis on the original closed pair. In this paper, we develop a closure comparison and a cut-and-reglue argument for this purpose.

Our main result is the following:

\begin{theorem}
Suppose that $K$ is a rationally null-homologous Floer simple knot in a closed irreducible 3-manifold with nonzero Thurston norm. Then $K$ is the unknot.
\end{theorem}

The proof begins by constructing balanced sutured manifolds from a taut surface in the knot complement. A closure comparison relates their sutured Floer homology to appropriate extremal Floer groups of a closed pair. We then transport the original Floer-simple equality through cutting and regluing, using generalized knot Floer homology to accommodate intermediate knots that are not rationally null-homologous.
\section{Preliminary}
All 3-manifolds and surfaces in this paper are assumed to be smooth and oriented. For an embedded submanifold $\Gamma$ of a manifold $N$,
we let $\nu \Gamma$ denote a closed regular neighborhood of $\Gamma$ in $N$.
For properly embedded submanifolds, neighborhoods are chosen
compatibly with the boundary. We write $N\setminus\!\setminus \Gamma$ as the closure of the complement of $\nu \Gamma$ in $N,$ with corners rounded.

\subsection{Sutured manifolds and sutured Floer homology}

A sutured manifold $(N,\gamma)$
consists of a compact oriented three-manifold together with
a possibly empty union of annuli $\gamma\subset\partial N$,
each equipped with an oriented core curve; the union of the core curves is
denoted by $s(\gamma)$ [3]. The sutures give a decomposition $\partial N=R_+(\gamma)\cup \gamma\cup R_-(\gamma),$ where $R_\pm(\gamma)$ are oriented coherently with the sutures,
with normal vectors pointing out of and into $N$, respectively.
The pair is \emph{balanced} if $N$ has no closed components,
every component of $\partial N$ meets $\gamma$, and $\chi(R_+(\gamma))=\chi(R_-(\gamma)).$
It is \emph{taut} if $N$ is irreducible and $R_\pm(\gamma)$ are
incompressible and Thurston-norm minimizing in $H_2(N,\gamma)$.
Here the Thurston norm is obtained by minimizing
$\chi_-(S)=\sum_i\max\{0,-\chi(S_i)\}$ over embedded oriented
representatives. Additionally, see Definitions 2.5~2.6 of [8] for the meaning of taut surfaces and decompositions.

For a closed connected 3-manifold $Y,$ and some integer $n>0,$ we may construct a balanced sutured manifold $Y(n)$ by removing $n$ pairwise disjoint 3-balls from $Y$ and then equipping each resulting sphere boundary with a suture. If $(N,\gamma)$ is a connected balaced sutured manifold then $N(n)=(N,\gamma)\# S^3(n)$ [3].

A balanced sutured diagram
$(\Sigma,\boldsymbol\alpha,\boldsymbol\beta)$ has no closed
surface components; there are equally many attaching curves in the two
collections, and each collection is separately linearly
independent in $H_1(\Sigma;\mathbb Q)$.
For an admissible such diagram, the complex $SFC(N,\gamma)$
is generated by
\[
\mathbb T_\alpha\cap\mathbb T_\beta
\subset\operatorname{Sym}^{g(\Sigma)}(\operatorname{int}\Sigma),
\qquad
\mathbb T_\alpha=\prod_{i=1}^{g(\Sigma)}\alpha_i,\quad
\mathbb T_\beta=\prod_{i=1}^{g(\Sigma)}\beta_i.
\]
Its differential counts Maslov-index-one holomorphic Whitney
disks modulo translation. Its homology splits over relative
$\operatorname{Spin}^{c}$ structures:
\[
SFH(N,\gamma)=
\bigoplus_{\xi\in\operatorname{Spin}^{c}(N,\gamma)}
SFH(N,\gamma,\xi).
\]
These structures are represented by nowhere-zero vector fields
with the prescribed sutured boundary behavior (see Section 4 of [3]).

\subsection{Domains and admissibility}

Let $(\Sigma,\boldsymbol\alpha,\boldsymbol\beta)$ be a
Heegaard diagram, and let $\Omega_1,\ldots,\Omega_n$
be the closures of the components of $\Sigma\setminus (\boldsymbol\alpha\cup \boldsymbol\beta)$. A \emph{domain} is an integral two-chain
\[
D=\sum_{j=1}^{n}m_j\Omega_j,\qquad m_j\in\mathbb Z.
\]
For $p$ in the interior of $\Omega_j$, its
\emph{local multiplicity} is $n_p(D)=m_j$.
The domain is nonnegative if all its coefficients are
nonnegative.

A \emph{periodic domain} is a domain $P$ whose boundary
is an integral combination of complete attaching circles:
\[
\partial P
=\sum_i m_{\alpha_i}\alpha_i+\sum_j m_{\beta_j}\beta_j.
\]
In a sutured diagram, the local multiplicities of $P$ vanish in every region meeting $\partial\Sigma$. In a pointed closed
diagram, we impose zero multiplicity at the designated
basepoints. Each attaching circle corresponds to a compressing disk in the ambient 3-manifold, and capping the boundary of $P$ with the corresponding
compressing disks, with the appropriate multiplicities
and orientations, gives a homology class $H(P)$ of the three-manifold represented by the diagram.

The diagram is \emph{admissible} if every nonzero
periodic domain has both positive and negative
coefficients.

For each $\Omega_j$, viewed as a surface
with $n_{\mathrm{c}}$ right-angled corners, its \emph{Euler measure} is $e(\Omega)=\chi(\Omega)-\frac{n_c}{4}$ (See [5]).
Extend this linearly to domains:
$e(D)=\sum_j m_j e(\Omega_j)$ for $D=\sum_j m_j\Omega_j$.

Let $\mathbf x,\mathbf y\in\mathbb T_\alpha\cap\mathbb T_\beta$
be generators. The point measure $n_{\mathbf x}(D)=\sum_{p\in\mathbf x}n_p(D)$; $n_{\mathbf y}(D)$ is computed similarly.

For a Whitney-disk class $\varphi$ from $\mathbf x$ to
$\mathbf y$ with domain $D(\varphi)$, its Maslov index is
\[
\mu(\varphi)=e(D(\varphi))+n_{\mathbf x}(D(\varphi))+n_{\mathbf y}(D(\varphi)).
\]

\section{Constructing Balanced Sutured Manifolds}
Suppose that $K$ is a rationally null-homologous Floer simple knot in a closed irreducible 3-manifold $Y$ with nonzero Thurston norm. First of all, we hope to transfer the problem about the closed pair $(Y,K)$ to a problem about some related balanced sutured manifold.

By discussion in [8], we may find a connected taut surface $G$ in $Y$ disjoint from $K$ with $g(G)>1.$ Let $M=Y\setminus\!\setminus G$, and denote
the two boundary copies of $G$ by $G_+$ and $G_-$, both inheriting the orientation of $G$, with normal vectors
pointing out of and into $M$, respectively. Set $E=M\setminus\!\setminus K$ and
$T_K=\partial\nu K$. Let $\rho\subset T_K$ consist of disjoint annular
neighborhoods of two oppositely oriented meridians of $K$.
Denote the positive and negative annuli of
$T_K\setminus\operatorname{int}\rho$ by $A_+$ and $A_-$,
oriented compatibly with the sutures. We equip $E$ with
the sutured structure
\[
R_+(\rho)=G_+\sqcup A_+,
\qquad
R_-(\rho)=G_-\sqcup A_-.
\]

Let $i_*:H_1(\partial E;\mathbb{Q})\to H_1(E;\mathbb{Q})$ be the map induced by inclusion. Since $\partial E\cong G_+\sqcup G_-\sqcup T_K,$ where $T_K$ is the boundary torus of a neighborhood of $K$, we have $H_1(\partial E;\mathbb{Q})=H_1(G_+;\mathbb{Q})\oplus H_1(G_-;\mathbb{Q})\oplus H_1(T_K;\mathbb{Q}).$ Let $W=H_1(G_+;\mathbb{Q})\oplus H_1(G_-;\mathbb{Q})\subset H_1(\partial M;\mathbb{Q}).$ We know that $\dim W=4g(G),$ $\dim H_1(\partial E;\mathbb{Q})=4g(G)+2$, and by Lemma 4.3 of [8], $\dim \ker(i_*)=2g(G)+1$. Since $g(G)>0$ we have $\dim W+\dim \ker(i_*)>\dim H_1 (\partial E;\mathbb{Q}).$ Hence, $W\cap \ker(i_*)$ is nontrivial. In fact, let $L_E=W\cap \ker(i_*)$, then $$\dim L_E\geq \dim W+\dim \ker(i_*)-\dim H_1 (\partial E;\mathbb{Q})=2g(G)-1.$$

Let $p_\pm:L_E\to H_1(G_\pm;\mathbb{Q})$ be the natural projections. Since $\ker p_+$ is an isotropic subspace of the symplectic vector space $H_1(G_+;\mathbb{Q})$, it has dimension at most $g(G)$. Similarly, $\ker p_-$ has dimension at most $g(G)$. By assumption, $2g(G)-1>g(G)$, so neither of $\ker p_\pm$ is all of $L_E.$ Moreover, a $\mathbb{Q}$-vector space is not the union of two proper subspaces, so we may choose $(b_+,-b_-)\in L_E$ with $b_\pm$ nontrivial. 

We next seek primitive classes
$c_\pm\in H_1(G_\pm;\mathbb Z)$ with equal images in
$H_1(E;\mathbb Q)$.
By Lemma 4.5 of [8], we may replace $E$
by a manifold with the same boundary and rational
homology, but torsion-free integral first homology.
The isomorphisms in Lemma 4.5 commute with boundary
inclusion, so the desired rational relation between
$c_+$ and $c_-$ is unaffected. Hence, from now on, we assume that $H_1(E;\mathbb{Z})$ is torsion-free; we retain this assumption until we finish proving the existence of such primitive classes.

Let $H_\pm=H_1(G_{\pm};\mathbb{Z})$, and let $V_\pm =\ker (H_\pm\to H_1(E;\mathbb{Z})).$ Let $K_E=\ker (H_1(\partial E;\mathbb{Z})\to H_1(E;\mathbb{Z}))$, and let $L_E^\mathbb{Z}=K_E\cap (H_+\oplus H_-\oplus 0)$. Let $P_\pm=\Pr_\pm(L_E^{\mathbb{Z}})\subset H_\pm.$ 
\begin{lemma}
There is a finite set of primes $\mathcal{P}=\{p_1,\cdots,p_n\}$ such that
\\1. If a prime $q\not\in \mathcal{P}$ and $c\in H_\pm$ maps to zero in $H_1(E;\mathbb{Z}/q\mathbb{Z}),$ then $c\in V_\pm +q H_\pm$.
\\2. If $c\in P_\pm$ is primitive in $P_\pm,$ and $c$ is not divisible by $p_i\in \mathcal{P}$ as an element of $H_\pm,$ then $c$ is primitive in $H_\pm.$
\end{lemma}
\begin{proof}
    Consider the long exact sequence of the pair $(E,G_+):$
    $$H_1(G_+)\to H_1(E)\to H_1(E,G_+)\to H_0(G_+)\to H_0(E).$$
    Since $G_+$ and $E$ are connected, the map $H_0(G_+)\to H_0(E)$ is an isomorphism. Therefore, $H_1(E)\to H_1(E,G_+)$ is surjective. Since $V_+=\ker (H_+\to H_1(E))$, we get the short exact sequence
    \begin{equation}0\to H_+/V_+\to H_1(E)\to H_1(E,G_+)\to 0.\end{equation}
    
Similar argument works for the ``-" sign.

    Now, choose $\mathcal{P}_1$ to contain all primes dividing the orders of $H_1(E,G_+)$, $H_1(E,G_-).$

    Let $q\not\in \mathcal{P}_1.$ Tensor (1) with $\mathbb{F}_q=\mathbb{Z}/q\mathbb{Z}$. Since $q$ avoids the torsion of $H_1(E,G_+)$, the Tor term vanishes, i.e. $\operatorname{Tor} (H_1(E,G_+),\mathbb{F}_q)=0.$ Hence,
    $$0\to (H_+/V_+)\otimes \mathbb{F}_q \to H_1(E)\otimes \mathbb{F}_q$$ is injective.

    Since $H_1(E)$ is torsion-free, $H_1(E)\otimes \mathbb{F}_q \cong H_1(E;\mathbb{F}_q).$ Also, $(H_+/V_+)\otimes \mathbb{F}_q\cong H_+/(V_++qH_+).$ Thus, we have an injection $H_+/(V_++qH_+)\hookrightarrow H_1(E;\mathbb{F}_q).$

    Therefore, the commutative diagram is:
    \[
\begin{tikzcd}[row sep=large, column sep=large]
H_\pm 
  \arrow[r, "i_\pm"] 
  \arrow[d] 
& H_1(E) 
  \arrow[d] 
\\
H_\pm/qH_\pm 
  \arrow[r, "i_{\pm,q}"] 
  \arrow[d] 
& H_1(E;\mathbb F_q) 
  \arrow[d, equal] 
\\
H_\pm/(V_\pm+qH_\pm) 
  \arrow[r, hook, "\bar i_{\pm,q}"] 
& H_1(E;\mathbb F_q).
\end{tikzcd}
\]
    Every arrow is induced either by inclusion $G_+\mapsto E,$ reduction mod $q,$ or quotienting by $V_++qH_+.$ The bottom map is well-defined because $V_+$ maps to zero in $H_1(E).$

    Now, suppose $c\in H_+$ maps to zero in $H_1(E;\mathbb{F}_q)$. Then its image in $H_+/(V_++qH_+)$ maps to zero under the injective map $H_+/(V_++qH_+)\hookrightarrow H_1(E;\mathbb{F}_q).$ Hence, $c\in V_++qH_+.$

    Now, we shall prove the second statement. Consider the short exact sequence
    \begin{equation}
        0\to P_+\to H_+\to H_+/P_+\to 0.
    \end{equation}
    Choose $\mathcal{P}_2$ to contain every prime dividing the orders of $H_+/P_+,$ $H_-/P_-.$ Let $\mathcal{P}=\mathcal{P}_1\cup \mathcal{P}_2.$ Now take $q\not\in \mathcal{P}.$ Tensor (2) with $\mathbb{F}_q$. Since $q$ avoids the torsion of $H_+/P_+,$ we have $\operatorname{Tor}(H_+/P_+,\mathbb{F}_q)=0.$ Therefore the induced map $P_+\otimes \mathbb{F}_q\to H_+\otimes \mathbb{F}_q$ is injective. Equivalently, the map $P_+/qP_+\hookrightarrow H_+/qH_+$ is injective.

    Now, suppose $c\in P_+$ is primitive in $P_+,$ and $c$ is not divisible in $H_+$ by any prime in $\mathcal{P}$. We want to show $c$ is primitive in $H_+.$

    Assume not. Then, there is a prime $q$ and an element $h\in H_+$ such that $c=qh.$ Since $c$ is not divisible by any bad prime, $q\not\in \mathcal{P}$. Also, $c$ is zero in $H_+/qH_+,$ but $c$ is primitive in $P_+$ so $c\not\in qP_+.$ Thus the class $c$ is nonzero in $P_+/qP_+,$ which contradicts that $P_+/qP_+\hookrightarrow H_+/qH_+$ is injective. Hence, $c$ is primitive in $H_+.$

    Similar argument works for the ``-" sign.
\end{proof}
\begin{lemma}
    Let $p$ be a prime number. Let $A_p^+=\{b=(b_+,b_-)\in L_E^\mathbb{Z}|b_+\in pH_+\}$, and let $\overline{A}_p^+$ be its image in $L_E^\mathbb{Z}/pL_E^\mathbb{Z}.$ Define $A_p^-$ and $\overline{A}_p^-$ similarly. Then, $\dim_{\mathbb{F}_p} \overline{A}_p^\pm\leq g(G).$
\end{lemma}
\begin{proof}
    Since $H_1(E)$ is torsion-free, $L_E^\mathbb{Z}$ is a direct summand of $W\subset H_1(\partial E).$ If the dimension of $\overline{A}_p^+$ is greater than $g(G),$ we can suppose that $(b_1^+,b_1^-),\cdots,(b_{g(G)+1}^+,b^-_{g(G)+1})$ are $(g(G)+1)$ elements, and they are linearly independent over $\mathbb{F}_p$ modulo $p H_1(\partial E).$ Let $V_{p,\pm}=\ker(H_1(G_\pm;\mathbb{F}_p)\to H_1(E;\mathbb{F}_p))$. As each $b_i^+$ is contained in $pH_+,$ each $b_i^-$ represents an element $[b_i^-]\in V_{p,-}.$ Also, the $[b_i^-]$s are independent in $V_{p,-}.$

    However, we can show that $V_{p,-}$ has dimension at most $g(G)$. Let $K_{E,p}=\ker (H_1(\partial E;\mathbb{F}_p)\to H_1(E;\mathbb{F}_p)).$ By Ni's Lemma 4.3 [8], over $\mathbb{F}_p,$ $K_{E,p}$ is Lagrangian in the symplectic vector space $H_1(\partial E;\mathbb{F}_p)$. Hence, $V_{p,-}$ is an isotropic subspace of $H_1(G_-;\mathbb{F}_p),$ which means $\dim_{\mathbb{F}_p} V_p^-\leq g(G).$

    All in all, we see that $\dim_{\mathbb{F}_p} \overline{A}_p^+\leq g(G).$ Similarly, $\dim_{\mathbb{F}_p} \overline{A}_p^-\leq g(G).$
\end{proof}

\begin{proposition}
    We may find primitive homology classes $c_\pm \in H_1(G_\pm)$ such that $c_+=c_-$ in $H_1(E;\mathbb{Q}).$
\end{proposition}
\begin{proof}
    Consider any $p\in \mathcal{P}$. Since $\dim_{\mathbb{F}_p} (L_E^\mathbb{Z}/pL_E^\mathbb{Z})= \operatorname{rank}L_E^\mathbb{Z}\geq 2g(G)-1>g(G)$, both $\overline{A}_p^+$ and $\overline{A}_p^-$ are proper subspaces. A vector space is not the union of two proper subspaces, so $(L^\mathbb{Z}_E/ pL_E^\mathbb{Z})\setminus (\overline{A}_p^+\cup \overline{A}_p^{-})\not=\emptyset$.

    Now, we consider three cases:

    \textbf{Case 1: $V_+\not=0$ and $V_-\not=0$}

    In this case, it suffices to choose primitive classes $c_+\in V_+$ and $c_-\in V_-.$

    \textbf{Case 2: $V_+=0$ and $V_-=0$}
    
    For each $p_i\in\mathcal{P},$ choose $r_i\in (L^\mathbb{Z}_E/ {p_i}L_E^\mathbb{Z})\setminus (\overline{A}_{p_i}^+\cup \overline{A}_{p_i}^{-})$. Use Chinese remainder theorem, we may choose $b=(b_+,b_-)\in L_E^\mathbb{Z}$ such that $b\equiv r_i$ mod $p_i L_E^\mathbb{Z}$ for each $p_i\in \mathcal{P}.$ Then neither $b_+$ nor $b_-$ is divisible by any bad prime $p_i.$

    Write $b_+=kc_+,$ $b_-=lc_-,$ where $c_\pm \in H_\pm$ are primitive and $k,l\in \mathbb{Z}_{\geq 0}.$ Let $d=\gcd(k,l),$ $k=dk',$ $l=dl',$ $\gcd(k',l')=1.$

    Since $L_E^\mathbb{Z}$ is saturated, we have $b'=(k'c_+,l'c_-)\in L_E^\mathbb{Z}$. Suppose $|l'|>1.$ Let $q$ be a prime divisor of $l'.$ Since $b'\in L_E,$ reducing mod $q$ gives $k'c_+\mapsto 0\in H_1(E;\mathbb{F}_q).$ Since $\gcd(k',l')=1,$ $q\nmid k'.$ Thus, $c_+\mapsto 0\in H_1(E;\mathbb{F}_q).$

    If $q\not\in \mathcal{P},$ then Lemma 3.1 and $V_+=0$ gives $c_+\in qH_+,$ contradicting the primitivity of $c_+.$

    Therefore, every prime divisor of $l'$ lies in $\mathcal{P}.$ However, that means $q|l,$ hence $q|b_-,$ contradicting the choice of $b.$ Thus, $|l'|=1,$ and similarly $|k'|=1.$ Therefore, $k=l=d.$ Dividing by $d$ over $\mathbb{Q}$, we get $c_+=c_-$ in $H_1(E;\mathbb{Q}).$

    \textbf{Case 3: The Mixed Case}

Suppose $V_+\not=0$ and $V_-=0.$

    For each prime $p,$ by Lemma 3.2, $(L_E^\mathbb{Z}/pL_E^\mathbb{Z})\setminus \overline{A}_p^-$ is nonempty. Then, again, applying Chinese remainder theorem, we could choose $b=(b_+,b_-)\in L_E^\mathbb{Z}$ such that $b_-\not=0,$ and $b_-\not\in p_i H_-$ for any $p_i\in \mathcal{P}.$

    Since $b_-\in P_-$ and $b_-\not=0,$ we may write $b_-=mc_-,$ where $c_-$ is primitive in $P_-.$

    Because $b_-$ is not divisible by any prime in $\mathcal{P},$ neither is $c_-.$ Hence, by Lemma 3.1, $c_-$ is primitive in $H_-.$

    Since $c_-\in P_-,$ there exists some $u_+\in H_+$ such that $(u_+,c_-)\in L_E^\mathbb{Z}.$ Also, because $V_+\not=0,$ and $V_+$ is a direct summand of $H_+,$ we may choose a primitive class $v_+\in V_+.$ Then $(v_+,0)\in L_E.$ By primitivity of $v_+,$ there is a homomorphism $\phi:H_+\to \mathbb{Z}$ with $\phi(v_+)=1.$ Take $n=1-\phi(u_+),$ and set $c_+=u_++nv_+.$ Then, $\phi(c_+)=1,$ and $(c_+,c_-)=(u_+,c_-)+n(v_+,0)\in L_E.$ Hence, $c_+,c_-$ are primitive and represent the same class in $H_1(E;\mathbb{Q}).$

    The case when $V_+=0$ and $V_-\not=0$ is analogous.
\end{proof}

Now, choose $c_\pm \in H_\pm$ as in Proposition 3.3. Let $k$ be the minimal positive integer such that $k[c_-]=k[c_+]$ in $H_1(E).$ Pick a curve $\lambda_+\in G_+$ representing $k[c_+]$ and a curve $\lambda_-\in G_-$ representing $k[c_-],$ such that $\lambda_+,\lambda_-$ are $k$ parallel copies of a non-separating simple closed curve. By assumption, the sutured manifold $(M,\emptyset)$ is taut, and hence $(E,\rho)$ is taut by [8]. By Theorem 2.8 of [8], we may choose a taut surface $S$ in $E\subset M$ such that $S$ has boundary $\lambda_+\sqcup(-\lambda_-)$ and $S$ has no closed components.

There is a homeomorphism $\phi:G_+\to G_-$ which sends $\lambda_+$ to $\lambda_-.$ Let $(Y',\overline{S})$ be the pair obtained from $(M,S)$ by gluing $G_+$ to $G_-$ via $\phi.$ Let $G^{(m)}$ be obtained from $\overline{S}$ and $m$ copies of $G$ by cut-and-paste. Let $K'$ be the new knot in $Y'.$ Then, $G^{(m)}$ is taut in $Y'\setminus\!\setminus K'$ since $S$ is taut in $E.$ By Lemma 5.6 of [8], $G^{(m)}$ is connected. Also, $G\cap \overline{S}$ provides product annuli connecting the two copies of $G^{(m)}$ in $M'=Y'\setminus\!\setminus G^{(m)}$. Thus $M'$ is connected.

Let $R_{\mathrm{out}}$ be the unique non-annular component of $G\setminus
\!\setminus(\overline{S}\cap G)$. Since $\overline{S}\cap G$ is the union of $k$ parallel copies of a
nonseparating curve, so $g(R_{\mathrm{out}})=g(G)-1\ge1$. Take \(m\ge k+1\). Denote the parallel copies of $G$ used in the construction of $G^{(m)}$ by $G_1,\cdots, G_m.$ Choose a product neighborhood
\(G\times I\subset Y'\setminus K'\) in which each $G_j=G\times \{t_j\}$, and suppose $t_1<\cdots<t_m.$ Let \(T_j\subset G^{(m)}\) be the
copy of \(R_{\mathrm{out}}\) inherited from \(G_j\). For \(1\le j\le k\), let \(Q_j\subset M'\) be the product region
corresponding to \(R_{\mathrm{out}}\times[t_j,t_{j+1}]\), with
horizontal boundary \(T_j\sqcup T_{j+1}\).
Splitting off \(Q_j\) from $M'$ along the two frontier annuli \(\partial R_{\mathrm{out}}\times[t_j,t_{j+1}]\) leaves
a sutured manifold \((N_j,\delta_j)\).

Suppose \(N_1,\ldots,N_k\) are all disconnected.
Since \(M'\) is connected, each $N_j$ has exactly two components,
one meeting each frontier annulus. Thus, for each $j,$
\(G^{(m)}\setminus\operatorname{int}T_j\) has two components, one bounded by each circle in $\partial T_j$.
Label \(\partial T_j=a_j\sqcup b_j\) consistently cross levels, and let \(U_j\) be the component of \(G^{(m)}\setminus\operatorname{int}T_j\) bounded by \(a_j\).

For each $j$, the component of \(N_j\) meeting the annulus joining \(a_j\)
to \(a_{j+1}\) has positive and negative boundary regions
\(U_j\) and \(U_{j+1}\). Also, by Lemma 2.11 of [8], \(N_j\setminus\operatorname{int}\nu K'\) is taut, so its positive and negative regions
are homologous norm minimizers; the additional boundary regions on \(\partial\nu K'\) are annuli
of Euler characteristic zero, so they do not affect \(\chi_-\).
Hence \(\chi_-(U_j)=\chi_-(U_{j+1})\). Since each \(U_j\) has one boundary component,
\(g(U_1)=\cdots=g(U_{k+1})\).

Order the components of \(G\setminus\!\setminus(\overline {S}\cap G)\)
cyclically as \(P_0=R_{\mathrm{out}},P_1,\ldots,P_{k-1}\),
so that smoothing joins \(P_r\) on \(G_j\) to \(P_{r+1}\)
on \(G_{j+1}\). Thus the chain from \(T_1\) passes through the annuli
\(P_1,\ldots,P_{k-1}\) on \(G_2,\ldots,G_k\), respectively,
before reaching \(T_{k+1}\).
These annuli and their smoothing connections form an annulus
\(A\) from \(b_1\) to \(a_{k+1}\), after relabeling if necessary (see Figure 2).
The connected union \(U_1\cup_{a_1}T_1\cup_{b_1}A\) avoids
\(\operatorname{int}T_{k+1}\) and has only the boundary
\(a_{k+1}\), so it equals \(U_{k+1}\). Consequently,
$g(U_{k+1})=g(U_1)+g(R_{\mathrm{out}})>g(U_1),$
a contradiction. Fix an index \(j\) for which \(N_j\) is
connected, and take \(Q_0=Q_j\) and
\((M_0,\gamma_0)=(N_j,\delta_j)\).

\begin{figure}[H]
    \centering
    \begin{minipage}[t]{0.48\linewidth}
    \captionsetup{font=scriptsize}
        \centering
        
        \includegraphics[width=0.9\linewidth]{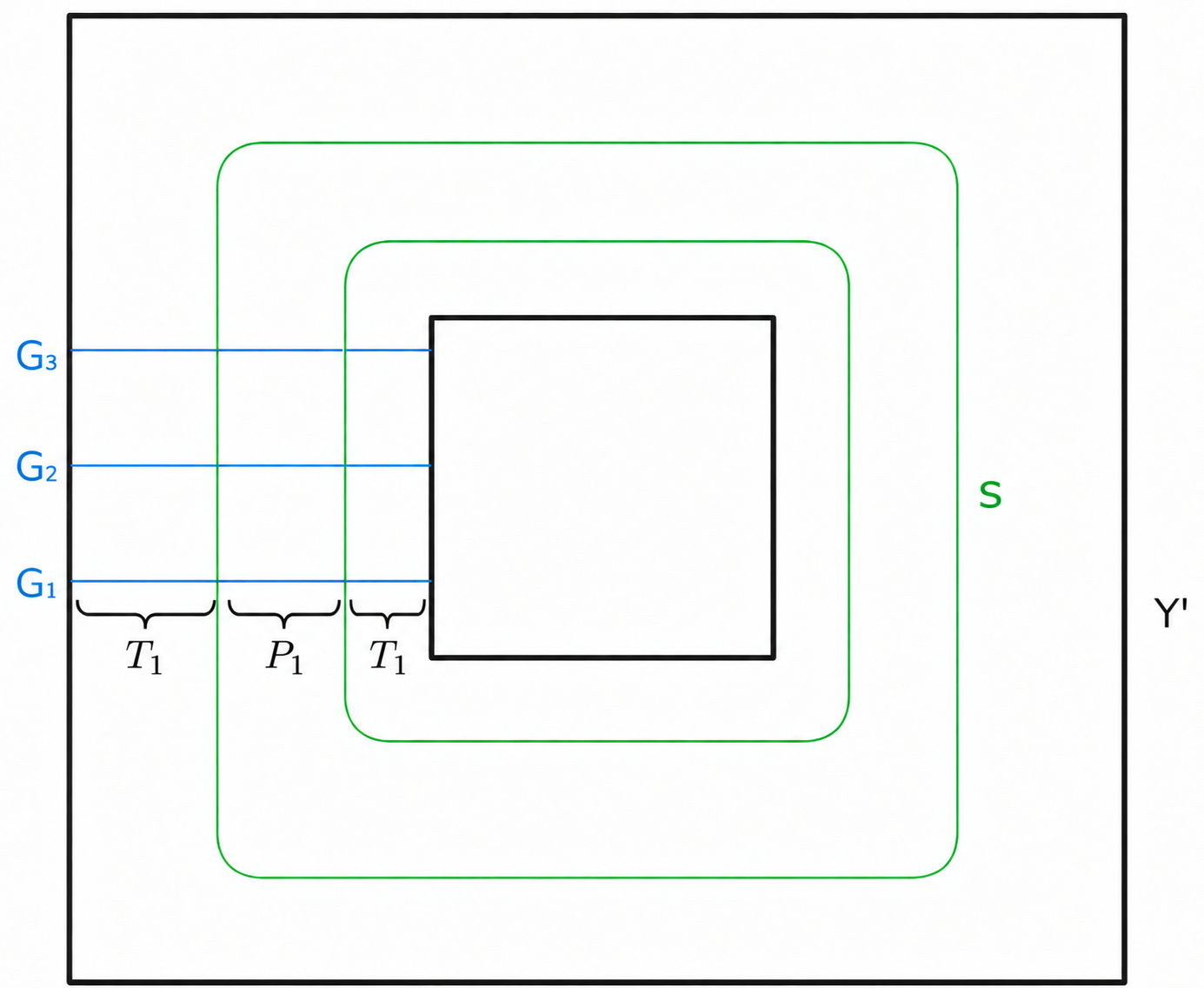}
        \caption{The configuration before smoothing: \(S\) and the parallel copies \(G_j\).}
        \label{fig:staircase}
    \end{minipage}\hfill%
    \begin{minipage}[t]{0.48\linewidth}
    \captionsetup{font=scriptsize}
        \centering
        
        \includegraphics[width=0.86\linewidth]{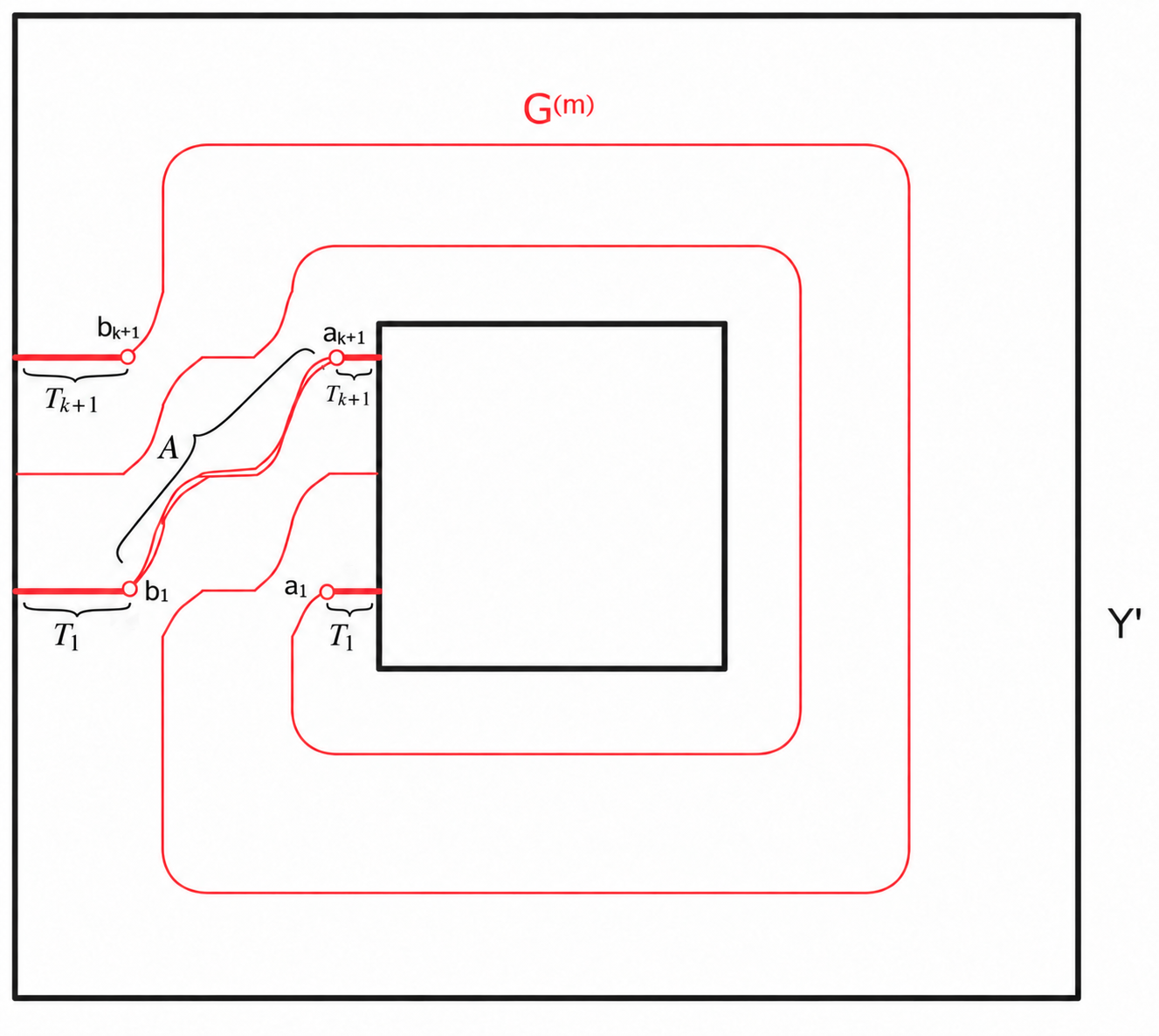}
        \caption{The smoothed surface \(G^{(m)}\), with the annulus \(A\subset G^{(m)}\) joining \(b_1\subset\partial T_1\) to \(a_{k+1}\subset\partial T_{k+1}\).
        }
        \label{fig:smoothing}
    \end{minipage}
\end{figure}

We also have the balanced sutured manifold $(M_0(K'),\gamma_0(K'))$, where $M_0(K')=M_0\setminus\!\setminus K'$ and $s(\gamma_0(K'))$ consists of $s(\gamma_0)$ and a pair of appropriately oriented meridians of $K'.$

By Lemma 2.11 of [8], $(M_0(K'),\gamma_0(K'))$ is taut because $G^{(m)}$ is taut in $Y'\setminus\!\setminus K'$, which makes the sutured manifold $(M'(K'),\mu(K'))$ taut. Here $M'(K')=M'\setminus\!\setminus K'$ and $s(\mu(K'))$ consists of a pair of appropriately oriented meridians of $K'$. Also, $M_0$ is irreducible. First, $Y'$ is irreducible. Given an embedded sphere $S^2$,
isotope $G\subset Y'$ to minimize the number of transverse intersection circles.
If any remain, an innermost intersection disk on $S^2$ is
boundary-parallel in $Y'\setminus\!\setminus G\cong M$, since $G_\pm$ are
incompressible and $M$ is irreducible. Isotoping $G$ across the
corresponding ball removes that intersection, contradicting minimality.
Thus $S^2$ is disjoint from $G$ and must bound a ball in $M$. Hence $Y'$ is irreducible. Now, let $U'=Y'\setminus M_0$. We know that $U'$ is connected because $\nu G^{(m)}$ meets $Q_0$. Now, let $S^2$ be an embedded sphere in $M_0.$ It must bound a 3-ball $B^3$ in $Y'.$ $U'$ is disjoint from $S^2$ so either $U'\subset \operatorname{int} B^3$ or $U'$ is disjoing from $B^3.$ If $U'\subset \operatorname{int} B^3$ then $[G^{(m)}]=0$ in $H_2(Y')$. Yet $[G^{(m)}]\not=[0]$ in $H_2(Y')$ because $M'$ is connected, and an arc in $M'$ joining a point on $G^{(m)}_+$ to its corresponding point on $G^{(m)}_-$ closes in $Y'$ to a curve intersecting $G^{(m)}$
transversely once. Hence $B^3$ is disjoint from $U'$, so
$B^3\subset M_0$. Thus $M_0$ is irreducible.

\section{Closure comparison}

Consider the balanced sutured manifolds $(M_0,\gamma_0),$ $(M_0(K'),\gamma_0(K'))$ we constructed in previous section. Now, we want to show that $\operatorname{rank }SFH(M_0(K'),\gamma_0(K'))=2\operatorname{rank }SFH(M_0,\gamma_0).$

The original Floer-simple hypothesis is formulated in terms of $\widehat{HF}$ and $\widehat{HFK}$ for closed manifolds and knots in them, but the groups above are invariants of balanced sutured manifold. We therefore seek to realize these sutured groups as extremal summands of Floer groups associated to a closed pair. The regluing construction provides a natural candidate: it produces the closed pair $(Y',K')$ along with a surface $G^{(m)}$ whose complement is related to $(M_0,\gamma_0)$ and $(M_0(K'),\gamma_0(K')).$ Although $(Y',K')$ is generally different from the original pair $(Y,K)$, it is obtained from $(Y,K)$ by cutting and regluing along $G.$ Hence, we may transport the original Floer-simple equality to $(Y',K')$ by proving some cut-and-reglue invariance of Floer groups.

Note that the reglued knot $K'\subset Y'$ need not remain rationally null-homologous. We therefore replace the Alexander-graded knot Floer homology used for the original rationally null-homologous knot $K\subset Y$ by the generalized doubly pointed knot Floer homology introduced by Sahamie [11].

\begin{definition}\label{def:hfk-hat}
Let $K\subset Y$ be a knot, and let $(\Sigma,\boldsymbol{\alpha},\boldsymbol{\beta})$
be a Heegaard diagram adapted to $K$: $K$ is isotopic to a curve on $\Sigma$
meeting $\beta_1$ transversely in one point and disjoint from the remaining
$\beta$-curves. Choose a small interval $I\subset K$ containing $K\cap\beta_1$ and
no other intersection with the attaching circles, and mark its endpoints by
basepoints $w,z$ with $\partial I=z-w$. The hat knot Floer complex
$\widehat{CFK}(Y,K)$ is the free $\mathbb{Z}_2$-module generated by the
intersection points $\mathbb{T}_\alpha\cap\mathbb{T}_\beta\subset
\mathrm{Sym}^{g(\Sigma)}(\Sigma)$, with differential
\[
\widehat{\partial}(\mathbf{x})
=\sum_{\mathbf{y}\in\mathbb{T}_\alpha\cap\mathbb{T}_\beta}
\sum_{\substack{\varphi\in\pi_2(\mathbf{x},\mathbf{y}),\ \mu(\varphi)=1\\
n_w(\varphi)=n_z(\varphi)=0}}
\#\widehat{\mathcal{M}}(\varphi)\;\mathbf{y},
\]
where $\#\widehat{\mathcal{M}}(\varphi)$ counts holomorphic disks in the class $\varphi$
modulo translation, and the diagram is chosen to be weakly admissible. Its homology
$\widehat{HFK}(Y,K)$ splits over the relative $\mathrm{Spin}^c$ structures
of the knot exterior and depends, up to canonical isomorphism, only on the pair
$(Y,K)$.
\end{definition}

In the language of Alishahi–Eftekhary, it is obtained from their refined sutured Floer complex for the meridionally sutured knot exterior by specializing the coefficient variables corresponding to the two meridional sutures to zero [1].

For a closed, connected, oriented surface
$G\subset Y\setminus K$ of genus greater than one, set
\[
B_Y([G]):=
\left\{
\mathfrak s\in\operatorname{Spin}^c(Y)
\;\middle|\;
\langle c_1(\mathfrak s),[G]\rangle=\chi(G)
\right\}.
\]
We call these the $G$-bottommost structures and write
\[
\begin{aligned}
\widehat{HF}(Y\mid G)
&:=
\bigoplus_{\mathfrak s\in B_Y([G])}
\widehat{HF}(Y,\mathfrak s),\\
\widehat{HFK}(Y,K\mid G)
&:=
\bigoplus_{\mathfrak s\in B_Y([G])}
\widehat{HFK}(Y,K,\mathfrak s).
\end{aligned}
\]
Here $\widehat{HFK}(Y,K,\mathfrak s)$ is the sum over
relative $\operatorname{Spin}^c$ structures of the knot
exterior that map to $\mathfrak s$ under meridional filling.
We use the analogous notation for $HF^+$ and for
$(Y',K')$ with either $G$ or $G^{(m)}$.

Now, we try to compare $SFH(M_0,\gamma_0)$ and $SFH(M_0(K'),\gamma_0(K'))$ with $ \widehat{HF}(Y'|G^{(m)})$ and $\widehat{HFK}(Y',K'|G^{(m)}).$

As stated in [7], being connected and
balanced, $(M_0,\gamma_0)$ admits a Morse function $h_0\colon M_0\to[0,1]$ with
only critical points of index $1$ and $2$ such that $h_0^{-1}(0)\cong R_-(\gamma_0)$ and
$h_0^{-1}(1)\cong R_+(\gamma_0)$. We can always arranged that $\Sigma_0:=h_0^{-1}(1/2)$ is the highest genus fibre which is connected and as one travels from $\frac{1}{2}$ to 0 one attaches 2-handles along $\beta_1,\cdots,\beta_d$ and as one travels from $\frac{1}{2}$ to 1 one attaches 2-handles along $\alpha_1,\cdots,\alpha_d,$ where the attaching circles are realized as vanishing circles of $h_0.$ The balanced condition means the set of $\alpha$ curves and respectively the set of $\beta$ curves are linearly independent in $H_1(\Sigma_0).$ Let $\boldsymbol{\alpha}_{\text{left}}=\{\alpha_1,\cdots,\alpha_d\}$ and $\boldsymbol{\beta}_{\text{left}}=\{\beta_1,\cdots,\beta_d\}$. The resulting tuple
$(\Sigma_0,\boldsymbol\alpha_{\text{left}},\boldsymbol{\beta}_{\text{left}})$ is a sutured Heegaard diagram for
$(M_0,\gamma_0)$ in the sense of Juh\'asz [4].

This diagram extends, in two steps, to a Heegaard diagram
of $(Y',f)$ with minimal fiber $\Sigma_{\min}=G^{(m)}$. First, complete
$(M_0,\gamma_0)$ back to $M'=Y'\setminus\!\setminus G^{(m)}$:
glue $T\times[0,1]$ to $M_0$ along
$\partial T\times[0,1]\cong A(\gamma_0)$, where $T\cong R_{\text{out}}$ defined in previous section, and extend $h_0$ over $T\times[0,1]$ as a product. At the level
of Heegaard diagrams it amounts to gluing the closing surface $T$ to
$\Sigma_0$ along their common boundary: the tuple
$(\Sigma_0\cup_\partial T,\boldsymbol\alpha_{\text{left}},\boldsymbol{\beta}_{\text{left}})$,
with the same attaching curves and no curve entering $T$, is a
Heegaard diagram of $M'$ viewed as a cobordism from $G^{(m)}_-$ to
$G^{(m)}_+$ --- indeed, compressing $\Sigma_0\cup T$ along the
$\alpha$-circles yields $R_+(\gamma_0)\cup_\partial T=G^{(m)}_+$, and
compressing along the $\beta$-circles yields $G^{(m)}_-$. Second,
identify the two boundary components by some homeomorphism, and the result
is $Y'$, and the Morse function descends to the circle-valued broken
fibration $f\colon Y'\to S^1$, whose fibers are connected because
$(M_0,\gamma_0)$ is balanced [7], and whose
extremal fibers are $\Sigma_{\min}=R_\pm\cup_\partial T=G^{(m)}$, of
genus $g$, and $\Sigma_{\max}=\Sigma_0\cup T$, of genus $g+d$. Lekili's
construction in Section 2.1 of [7] applied to $(Y',f)$ now produces the
Heegaard diagram used below: construct $\Sigma_*$ by tubing $\Sigma_{\max}=\Sigma\cup T$ to
$\Sigma_{\min}=G^{(m)}$ along two gradient flow lines, one over each
semicircle, and add additional attaching curves $\boldsymbol{\alpha}_{\text{right}},\boldsymbol{\beta}_{\text{right}}$ to $\boldsymbol{\alpha}_{\text{left}},\boldsymbol{\beta}_{\text{left}}$, where $\boldsymbol{\alpha}_{\text{right}},\boldsymbol{\beta}_{\text{right}}$ consist of the tube-seam curves $\alpha_0,\beta_0$ and other right-side curves which come from a dual arc system on the $\Sigma_{\min}$-side completed by the flow. Let $\boldsymbol{\alpha}_*=\boldsymbol{\alpha_{\text{left}}}\cup \boldsymbol{\alpha_{\text{right}}}$ and $\boldsymbol{\beta}_*=\boldsymbol{\beta_{\text{left}}}\cup \boldsymbol{\beta_{\text{right}}}$. The diagram $(\Sigma_*,\boldsymbol{\alpha}_*,\boldsymbol{\beta}_*)$ corresponds to the broken fibration $(Y',f).$ Let $F$ be the part of $\Sigma_*$ which consists of $\Sigma_{\max}$ (with two disks removed) together with halves of the connecting tubes up to $\alpha_0$ and $\beta_0.$ Let $\overline{F}=\Sigma_*\setminus \text{Int}(F).$ After winding along some closed curves in $F$ if necessary, the diagram can be made weakly admissible; see Lemma 8 of [7].

By Theorem 24 of [7], $SFH(M_0,\gamma_0)\cong \oplus_{\mathfrak{s}\in \mathcal{B}_{Y'}([G^{(m)}])} HF^+(Y',\mathfrak{s})$.

For the knot version, since the omitted product
piece is chosen disjoint from \(K'\), we have
\(K'\subset\operatorname{Int}(M_0)\). Via some isotopy supported in $M_0,$ we can arrange that $K'$ is in finite bridge position. We may perform a finite number of stabilizations, supported in
\(\operatorname{Int}(\Sigma_0)\). Each move is a meridional
stabilization of the ambient Heegaard splitting and decreases the
bridge number by one. Ultimately, $(\Sigma_*,\boldsymbol{\alpha}_*,\boldsymbol{\beta}_*)$ becomes $({\Sigma}'_*,\boldsymbol{\alpha}_*\cup \boldsymbol\alpha'_*,\boldsymbol{\beta}_*\cup \boldsymbol{\beta}'_*)$ after the stabilizations, and $K'$ becomes 1-bridge with respect to the corresponding Heegaard splitting of $Y'.$ That is, we can pick two basepoints $w_k,z_k$ in $\Sigma_0\setminus (\boldsymbol{\alpha}_*\cup \boldsymbol{\beta}_*)$. Connect $w_k$ to $z_k$ by an arc $a$ in $\Sigma_0\setminus \boldsymbol{\alpha}_*$, and push the interior of the resulting arc into the $\alpha$-compression bodies; similarly, connect $z_k$ to $w_k$ by another arc $b$ in $\Sigma_0\setminus \boldsymbol{\beta}_*$, and push the interior of the resulting arc into the $\beta$-compression bodies. We may arrange so that $K'=a\cup b$.

For later comparison with
Sahamie's doubly pointed complex, we perform one more stabilization
near \(w_k,z_k\): add a 1-handle with one foot near $w_k$ and one foot near $z_k,$ let
\(\alpha_s\) run along \(a\) and through the new handle, and let
\(\beta_s\) be the corresponding belt circle.  After isotopy, \(K'\)
lies on the stabilized surface, meets \(\beta_s\) once, and is
disjoint from all the other \(\beta\)-curves.

Now, after all those stabilizations that reduce bridge number and the stabilization introducing $\alpha_s,\beta_s$, suppose $\Sigma_*$ becomes $\widetilde{\Sigma}_*$, $\Sigma_0$ becomes $\widetilde{\Sigma}_0,$ $\Sigma_{\max}$ becomes $\widetilde{\Sigma}_{\max},$ and $F$ becomes $\widetilde{F}.$ Suppose the genus of $\widetilde{F}$ is $g+d'.$

Set
\[
\boldsymbol\alpha^s_\text{left}
=\boldsymbol\alpha_{\mathrm{left}}\cup \boldsymbol{\alpha}'_*\cup\{\alpha_s\},
\qquad
\boldsymbol\beta^s_{\text{left}}
=\boldsymbol\beta_{\mathrm{left}}\cup \boldsymbol{\beta}'_*\cup\{\beta_s\},
\] and put
\[
\Sigma_{K'}
=
\widetilde\Sigma_0
\setminus\operatorname{int}(D_{w_k}\cup D_{z_k}),
\]
where $D_{w_k}, D_{z_k}$ are small disks on $\widetilde{\Sigma}_0$ near $w_k,z_k,$ respectively. Let
\[
M_0=C_\alpha\cup_{\widetilde\Sigma_0}C_\beta
\]
be the compression-body splitting corresponding to the tuple $(\widetilde{\Sigma}_0,\boldsymbol{\alpha}^s_{\text{left}},\boldsymbol{\beta}^s_{\text{left}}).$

We know that $K'=a\cup b$ with $a,b$ boundary parallel. Thus, there are arcs $c_a,c_b$ on $\widetilde{\Sigma}_0$ with endpoint $z_k,w_k$ such that $a$ and $c_a$ cobound a bridge disk $\Delta^r_\alpha$, $b$ and $c_b$ cobound a bridge disk $\Delta^r_\beta$. Let $B^r_\alpha$ be a regular neighborhood of $\Delta^r_\alpha$ in $C_\alpha$ and $B^r_\beta$ be a regular neighborhood of  $\Delta^r_\beta$ in $C_\beta$, with $B^r_\alpha\cap \widetilde{\Sigma}_0=E^r_\alpha$ and $B^r_\beta\cap \widetilde{\Sigma}_0=E^r_\beta$. Then, $E^r_\alpha$ (resp. $E^r_\beta$) is a small disks containing $c_a$ (resp. $c_b$) together with
\(D_{w_k}\cup D_{z_k}\).  

Set $\alpha_r=\partial E^r_\alpha$ and $\beta_r=\partial E^r_\beta$ (see Figure 3). Although \(\alpha_r,\beta_r\) both bound a disk in $\widetilde{\Sigma}_0$, they cut
off the pairs of pants
\[
E^r_\alpha\setminus
\operatorname{int}(D_{w_k}\cup D_{z_k})
\] and \[
E^r_\beta\setminus
\operatorname{int}(D_{w_k}\cup D_{z_k})
\]
on \(\Sigma_{K'}\), respectively.
\begin{figure}[H]
 \centering
  \includegraphics[width=0.5\textwidth]{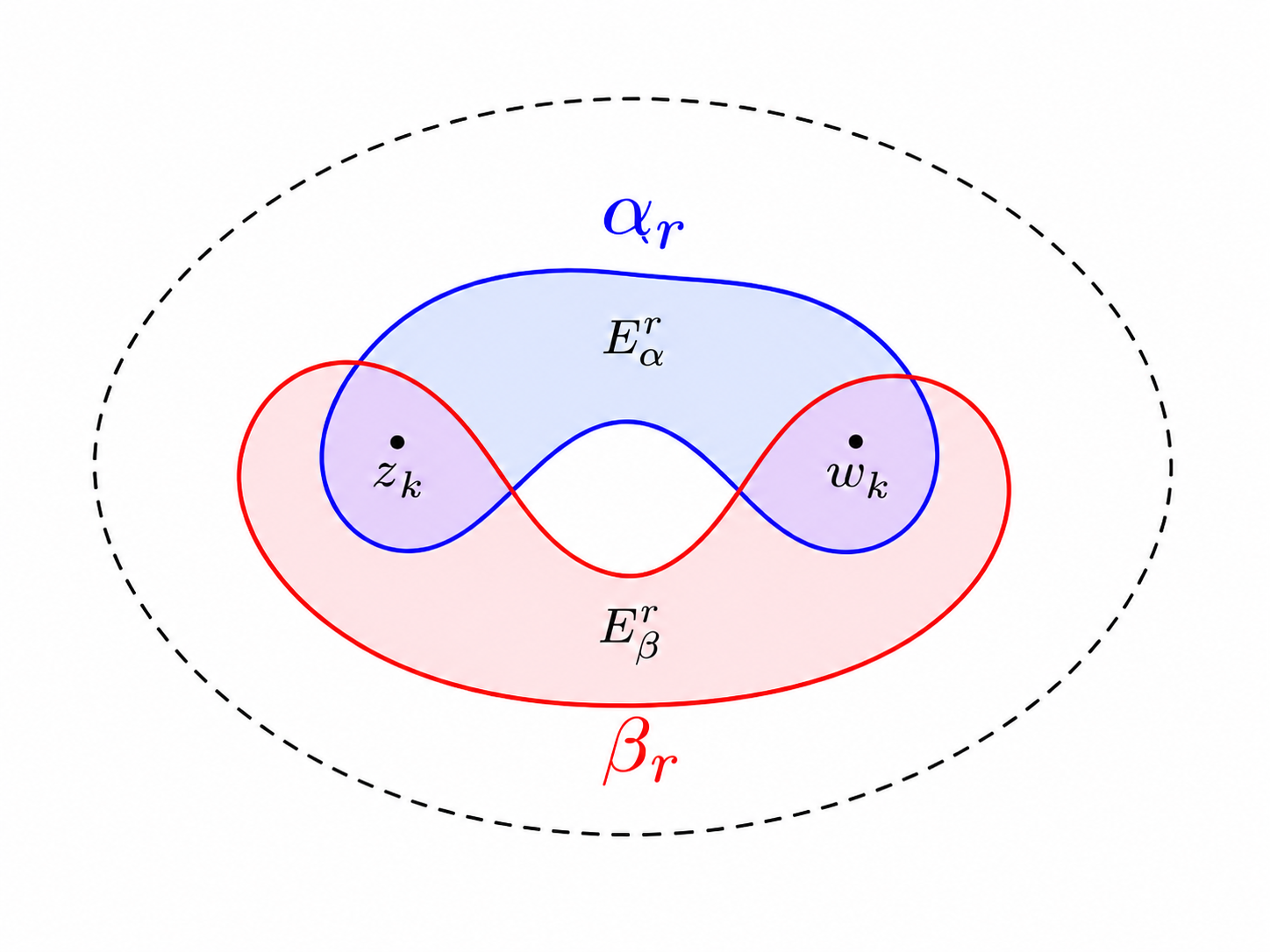}
 \caption{The local configuration of \(\alpha_r\) and \(\beta_r\).}
 \end{figure}

Consider $D^r_\alpha=\partial B^r_\alpha\setminus \text{int}(E^r_\alpha)$ and  $D^r_\beta=\partial B^r_\beta\setminus \text{int}(E^r_\beta)$. Push $D^r_\alpha,D^r_\beta$ slightly to the interior, and they become properly embedded disks in \(C_\alpha\setminus\!\setminus a\) and \(C_\beta\setminus\!\setminus b\), with boundaries $\alpha_r,\beta_r,$ respectively.

Compressing along  \(D^r_\alpha\) (resp. $D^r_\beta$) turn $\Sigma_{K'}$ into two disjoint components. One component is obtained by capping the pair of pants $
E^r_\alpha\setminus
\operatorname{int}(D_{w_k}\cup D_{z_k})$ (resp. $
E^r_\beta\setminus
\operatorname{int}(D_{w_k}\cup D_{z_k})$) into an
annulus parallel to
\(\partial\nu(a)\cap C_\alpha\) (resp. \(\partial\nu(b)\cap C_\beta\)). The other component is just homeomorphic to the original $\widetilde{\Sigma}_0$. 
Consequently, the old disk systems together with
\(D^r_\alpha,D^r_\beta\) are complete disk systems for the two
compression bodies obtained by drilling \(a,b\).
The new boundary components
\(\partial D_{w_k},\partial D_{z_k}\) determine the two meridional
suture annuli.  Hence
\[
\left(
\Sigma_{K'},\boldsymbol{\alpha}^s_{\text{left}}\cup \{\alpha_r\},\boldsymbol{\beta}^s_{\text{left}}\cup \{\beta_r\}
\right)
\]
is a balanced sutured Heegaard diagram for
\[
\bigl(M_0(K'),\gamma_0(K')\bigr).
\]

Complete
$(M_0(K'),\gamma_0(K'))$ back to $M'\setminus\!\setminus K'$:
glue $T\times[0,1]$ to $M_0(K')$ along
$\partial T\times[0,1]\cong A(\gamma_0)$, where $T\cong R_{\text{out}}$ as before, and extend the restriction of $h_0$ to $M_0(K')$ over $T\times[0,1]$ as a product. At the level
of Heegaard diagrams it amounts to gluing the closing surface $T$ to
$\Sigma_{K'}$ along their common boundary: the tuple
$({\Sigma}_{K'}\cup_\partial T,\boldsymbol{\alpha}^s_{\text{left}}\cup \{\alpha_r\},\boldsymbol{\beta}^s_{\text{left}}\cup \{\beta_r\})$,
with the same attaching curves and no curve entering $T$, is a
sutured diagram of $M'\setminus\!\setminus K'$ --- indeed, compressing $\Sigma_{K'}\cup T$ along the
$\alpha$-circles yields $R_+(\gamma_0(K'))\cup_\partial T=G^{(m)}_+\cup \partial\nu(c_a)$, and
compressing along the $\beta$-circles yields $G^{(m)}_-\cup\partial\nu(c_b)$. 

Now, construct a surface $\widetilde{\Sigma}_{K'}$ by tubing \(\Sigma_{K'}\) to \(G^{(m)}\) along two gradient flow lines, and adjoin the boundaries of the
tube-meridian disks and the disks swept out by standard arc systems
on \(G^{(m)}\). These are precisely
\(\boldsymbol\alpha_{\rm right}\) and
\(\boldsymbol\beta_{\rm right}\). Decomposing $\widetilde{\Sigma}_{K'}$ along $\boldsymbol{\alpha}^s_{\text{left}}\cup \{\alpha_r\}\cup\boldsymbol{\alpha}_{\text{right}}$ gives us $\partial\nu(c_a)\sqcup S^2.$ Decomposing $\widetilde{\Sigma}_{K'}$ along $\boldsymbol{\beta}^s_{\text{left}}\cup \{\beta_r\}\cup\boldsymbol{\beta}_{\text{right}}$ gives us $\partial\nu(c_b)\sqcup S^2.$ 
Choose a free basepoint $w_r$ on $\Sigma_{\min}$. We may assume that some section of $f$ intersects $\Sigma_{\max}$ at $w_k$ and intersects $\Sigma_{\min}$ at $w_r$. Deleting a disk $D_{w_r}$ about the free basepoint \(w_r\) turns these spherical
components into disks \(D_+\) and \(D_-\).  Hence the final compressed
surfaces are
\[
\partial\nu(c_a)\sqcup D_+,
\qquad
\partial\nu(c_b)\sqcup D_-,
\]
while the three components of the vertical boundary carry the two
meridional sutures on \(\partial\nu K'\) and the single suture $\tau$ on the
new spherical boundary. Therefore, let $\widetilde{\Sigma}_{K'}^\circ=\widetilde{\Sigma}_{K'}\setminus D_{w_r};$ let $X^\circ=Y'\setminus\operatorname{int}(\nu K'\cup B^3);$ let $\boldsymbol{\alpha}^{K'}_{\text{left}}=\boldsymbol{\alpha}^s_{\text{left}}\cup \{\alpha_r\}$ and $\boldsymbol{\beta}^{K'}_{\text{left}}=\boldsymbol{\beta}^s_{\text{left}}\cup \{\beta_r\}$. Then the completed diagram 

\[
\widetilde{\mathcal{H}}^\circ_{K'}=\left(
\widetilde{\Sigma}^\circ_{K'},\boldsymbol{\alpha}^{K'}_{\text{left}}\cup\boldsymbol{\alpha}_{\text{right}},\boldsymbol{\beta}^{K'}_{\text{left}}\cup\boldsymbol{\beta}_{\text{right}}
\right).
\]

represents
\[
\left(
X^\circ,
\mu_{K'}\cup\tau
\right).
\]

If $\widetilde{\mathcal{H}}^\circ_{K'}$ is weakly admissible, then it computes $SFH(X^\circ,
\mu_{K'}\cup\tau)$. Moreover,
\[
SFH\!\left(
X^\circ,
\mu_{K'}\cup\tau
\right)
\cong
\widehat{HFK}(Y',K')\otimes V.
\]
Indeed, removing the two knot basepoints from a Sahamie-adapted
doubly pointed diagram identifies his generalized knot Floer complex
with the sutured Floer complex of the meridionally sutured knot exterior,
whereas removing the additional ball contributes the standard
rank-two factor \(V\).  If the ball is chosen disjoint from
\(G^{(m)}\), this isomorphism preserves the \(G^{(m)}\)-degree.

Denote $\mu_{K'}\cup \tau$ by $\gamma^\circ$. For a relative class $\xi\in\text{Spin}^c(X^\circ,\gamma^\circ)$,
write $\xi\in\text{Spin}^c(X^\circ,\gamma^\circ\mid G^{(m)})$ when $\langle c_1(\xi),[G^{(m)}]\rangle=\chi(G^{(m)}).$
Since \(G^{(m)}\) is closed and contained in the interior, the displayed
evaluation is independent of any boundary trivialization.

We say that \(\widetilde{\mathcal{H}}^\circ_{K'}\) is \(G^{(m)}\)-admissible if it is weakly
\(\xi\)-admissible for every
\(\xi\in\text{Spin}^c(X^\circ,\gamma^\circ\mid G^{(m)})\); namely, every
nonzero periodic domain \(Q\) satisfying
\[
\left\langle c_1(\xi),H(Q)\right\rangle=0
\]
has both positive and negative coefficients.

\begin{proposition}
After winding along some closed curves in $\widetilde{F}$ if necessary, $\widetilde{\mathcal{H}}^\circ_{K'}$ is \(G^{(m)}\)-admissible.
\end{proposition}

\begin{proof}
By Lemma 7 of [7], we have a linearly independent collection $\mathcal{A}=\{[\alpha_1],\cdots,[\alpha_d],[\xi_1],\cdots,[\xi_{2g}],\delta\}$ of classes in $H_1(F,\partial F),$ where $\alpha_1,\cdots,\alpha_d$ are curves in $\boldsymbol{\alpha}_{\text{left}}$, the arcs $\xi_1,\cdots,\xi_{2g}$ are the $F$-portions of curves in $\boldsymbol{\alpha}_{\text{right}}$, and $\delta\subset F$ is a properly embedded arc joining
\(\alpha_0\) to \(\beta_0\), disjoint from all other $\alpha$ and $\beta$ curves. When taking complements, we identify collections of homology classes with the unions of their chosen curve and arc representatives.

We know that $\widetilde{F}$ is obtained from $F$ by iteratively doing connected sum operation with a torus. It is easy to arrange that the stabilizations take place away from $\delta$. We may also assume that the meridional stabilizations are performed near the shadow of some $\alpha$-bridge of $K'$, and thus we may arrange that the connected-sum disks $D_1,\cdots,D_q$ are all on $F,$ and $\widetilde{F}$ is constructed from $F$ by replacing each $D_i$ with a punctured torus $T_i.$ Each $T_i$ corresponds to the tube added after stabilization; let $\alpha_{T_i}$ be the new $\alpha$ curve introduced by the corresponding stabilization. Each $\alpha_{T_i}$ must be essential and hence nonseparating, so each $T_i\setminus \alpha_{T_i}$ is connected. Moreover, by Lekili's argument in [7], $F\setminus (\mathcal{A}\cup D_1\cup \cdots \cup D_q)$ is connected. Therefore, $\widetilde{F}\setminus(\mathcal{A}\cup\{[\alpha]:\alpha \in \boldsymbol{\alpha}'_*\cup \{\alpha_s\}\})$ is connected, and hence
\[\mathcal{A}'=\mathcal{A}\cup\{[\alpha]:\alpha \in \boldsymbol{\alpha}'_*\cup \{\alpha_s\}\}\] is a linearly independent collection of classes in $H_1(\widetilde{F},\partial\widetilde{F}).$ 

Now, let $\mathcal{H}_*=(\Sigma_*,\boldsymbol{\alpha}_*,\boldsymbol{\beta}_*,w_r
)$ and $\mathcal{H}^s_*=(\widetilde{\Sigma}_*,\boldsymbol{\alpha}^s_{\text{left}}\cup \boldsymbol{\alpha}_{\text{right}},\boldsymbol{\beta}_{\text{left}}^s\cup \boldsymbol{\beta}_{\text{right}},w_r
)$. By the argument above, for each $\alpha_t\in \boldsymbol{\alpha}^s_{\text{left}}\cup \boldsymbol{\alpha}_{\text{right}},$ Poincaré--Lefschetz duality supplies a closed curve $\zeta_{\alpha_t}\subset \widetilde{F}\setminus \delta$ such that $\zeta_{\alpha_t}$ has algebraic intersection number 1 with $\alpha_t$ and 0 with all other $\alpha$-curves of $\mathcal{H}^s_*.$ Then, we can apply similar winding argument as in Lemma 8 of [7]. We can make $\mathcal{H}^s_*$ weakly admissible for $\mathfrak{s}\in \mathcal{B}_{Y'}([G^{(m)}])$ by winding successively along the curves $\{\zeta_{\alpha_t}\}$ in $\widetilde{F}$. At each step, wind a distinct $\alpha_t\in \boldsymbol{\alpha}^s_{\text{left}}\cup \boldsymbol{\alpha}_{\text{right}}$ along its corresponding $\zeta_{\alpha_t}.$ Each winding is a diffeomorphism of $\widetilde{F}$ which is isotopic to the identity and supported near $\zeta_{\alpha_t}$, so the winding operation preserves disjointness of the $\alpha$-curves and the algebraic intersection pattern of the $\alpha$-curves with the $\zeta$-curves. Also, the winding operation does not change the bridge number of $K'$, so we can first make $\mathcal{H}^s_*$ weakly admissible for $\mathfrak{s}\in \mathcal{B}_{Y'}([G^{(m)}])$ after winding sufficiently many times, and then choose $z_k,w_k,\alpha_r,\beta_r$ as in preceeding discussion.

It remains to analyze the effect of adjoining
\(\alpha_r,\beta_r\).

Let $\{\Omega_i\}$ be the collection of components of $\Sigma^\circ_{K'}\setminus(\boldsymbol{\alpha}^{K'}_{\text{left}}\cup\boldsymbol{\alpha}_{\text{right}}\cup\boldsymbol{\beta}^{K'}_{\text{left}}\cup\boldsymbol{\beta}_{\text{right}})$. If \(Q\) is a periodic domain of $\widetilde{{H}}^\circ_{K'}$, then $Q$ is an integral combination of the regions in $\{\Omega_i\}$ and its boundary is an integral combination of the attaching circles and
has no component on \(\partial\widetilde{\Sigma}^\circ_{K'}=\partial D_{z_k}\cup \partial D_{w_k}\cup \partial D_{w_r}\).  Since the attaching
circles are contained in the interior of \(\widetilde{\Sigma}^\circ_{K'}\), each
component \(\partial D_p\), \(p\in\{w_k,z_k,w_r\}\), has a collar neighborhood contained in a region $\Omega_i\subset \widetilde{\Sigma}^\circ_{K'}\setminus Q$.  The coefficient of that $\Omega_i$ in $Q$ must be zero.  Consequently, \(Q\) extends uniquely
over the three capping disks with coefficient zero.  Denoting the
resulting two-chain on $\widetilde{\Sigma}_*$ by \(\overline Q\), we have
\[
\partial\overline Q=\partial Q,
\qquad
n_{w_k}(\overline Q)
=
n_{z_k}(\overline Q)
=
n_{w_r}(\overline Q)
=0.
\]

Note that the two disks $E^r_\alpha, E^r_\beta$ are on $\widetilde{\Sigma}_*$ and have boundary $\alpha_r,\beta_r$. Also, $n_{w_r}(E^r_\alpha)=n_{w_r}(E^r_\beta)=0.$ If the
coefficients of \(\alpha_r\) and \(\beta_r\) in \(\partial \overline{Q}\) are
\(m_1\) and \(m_2\), respectively, set
\[
P=\overline{Q}-m_1E^r_\alpha-m_2E^r_\beta.
\]
Then \(P\) is a periodic domain for $\mathcal{H}_*^s,$ because its boundary is an integral combination of attaching curves of $\mathcal{H}_*^s$, and $n_{w_r}(P)=0.$

If \(P=0\), then
\[
\overline{Q}=m_1E^r_\alpha+m_2E^r_\beta.
\]
Since \(\overline{Q}\) has zero multiplicity at \(w_k\) and
\(w_k\in E^r_\alpha\cap E^r_\beta\), we have \(m_1+m_2=0\).  Hence
\[
\overline{Q}=m_1(E^r_\alpha-E^r_\beta),
\]
which has both positive and negative coefficients unless \(m_1=0\).

Suppose now that \(P\neq0\) and $\langle c_1(\xi),H(Q)\rangle=0.$

Let \(j:X^\circ\hookrightarrow Y'\) be the natural inclusion. Write
\[
\partial\overline Q=m_1\alpha_r+m_2\beta_r+\eta,
\]
where \(\eta\) is supported on the old attaching curves in $\mathcal{H}_*^s$. Then \(\partial P=\eta\).

Suppose \[\eta=\sum_{\alpha_t\in \boldsymbol{\alpha}^s_{\text{left}}\cup \boldsymbol{\alpha}_{\text{right}}}u_t \alpha_t+\sum_{\alpha_t\in \boldsymbol{\beta}^s_{\text{left}}\cup \boldsymbol{\beta}_{\text{right}}}v_t \beta_t.\]  Let $D_\alpha^t$ denote the core of the two-handle
attached to $\alpha_t$ and $D_\beta^t$ the core of the two-handle attached to $\beta_t$. Let 
\[P_\eta=\sum_{\alpha_t\in \boldsymbol{\alpha}^s_{\text{left}}\cup \boldsymbol{\alpha}_{\text{right}}}u_t D_\alpha^t+\sum_{\alpha_t\in \boldsymbol{\beta}^s_{\text{left}}\cup \boldsymbol{\beta}_{\text{right}}}v_t D_\beta^t.\]

We know that the class $H(Q)$ is represented by $Q+P_\eta+m_1D_\alpha^r+m_2 D_\beta^r,$ and $j_*(H(Q))$ is represented by $\overline{Q}+P_\eta+m_1D_\alpha^r+m_2 D_\beta^r.$ Also, the class $H(P)$ is represented by $P+P_\eta.$ Hence, $j_*(H(Q))-H(P)$ is represented by
\begin{align*}
    \overline{Q}+m_1D_\alpha^r+m_2D_\beta^r-P&=m_1 E_\alpha^r+m_2 E_\beta^r+m_1D_\alpha^r+m_2D_\beta^r
    \\&=m_1\partial B_\alpha^r+m_2\partial B_\beta^r
    \\&=\partial (m_1 B_\alpha^r+m_2 B_\beta^r).
\end{align*}

Therefore
\[
j_*(H(Q))=H(P)\in H_2(Y').
\]
If \(\mathfrak s\) is the ambient extension of
\(\xi\in\text{Spin}^c(X^\circ,\gamma^\circ)\), then
\[
\langle c_1(\xi),H(Q)\rangle
=
\langle c_1(\mathfrak s),H(P)\rangle,
\]

and thus
\[
\langle c_1(\mathfrak{s}),H(P)\rangle=0.
\]

 Then, since we have arranged that $\mathcal{H}_*^s$ is weakly admissible for $\mathfrak{s}\in \mathcal{B}_{Y'}([G^{(m)}])$, the 2-chain $P$ has both positive and negative coefficient. Therefore, $\overline{Q}$ and thus \(Q\) also has both positive and negative coefficients.
This proves $G^{(m)}$-admissibility of $\widetilde{H}^\circ_{K'}$.
\end{proof}

 Hence, although $\widetilde{H}^\circ_{K'}$ might not let us compute the whole $SFH(X^\circ,\gamma^\circ),$ we may define
\[
SFH(X^\circ,\gamma^\circ\mid G^{(m)})
=
\bigoplus_{\substack{
\xi\in\text{Spin}^c(X^\circ,\gamma^\circ)\\
\langle c_1(\xi,t_{\gamma^\circ}),[G^{(m)}]\rangle
=\chi(G^{(m)})}}
SFH(X^\circ,\gamma^\circ,\xi).
\]
Since \(G^{(m)}\) is a closed surface contained in the interior of
\(X^\circ\), its Chern evaluation is independent of the choice of
boundary trivialization. By Proposition 4.2, we may compute $SFH(X^\circ,\gamma^\circ\mid G^{(m)})$ from $\widetilde{H}_{K'}^\circ$.

Recall that \[
\widehat{HFK}(Y',K'\mid G^{(m)})
=
\bigoplus_{\mathfrak{s}\in \mathcal{B}_{Y'}(G^{(m)})}
\widehat{HFK}(Y',K',\mathfrak s).
\]

Let $X=Y'\setminus\!\setminus K'$. We know that removing neighborhoods of \(w_k,z_k\) from the Heegaard surface gives the tautological
chain-level identification
\[
SFC(X,\mu_{K'})\cong\widehat{CFK}(Y',K').
\]
Since \(G^{(m)}\) is disjoint from the knot-filling region, this
identification preserves the Chern evaluation on \([G^{(m)}]\).

Now, consider $SFH(X^\circ,\gamma^\circ)$. Juhász's puncturing formula in [3] gives
\[
SFH(X^\circ,\gamma^\circ)
\cong SFH(X,\mu_{K'})\otimes V.
\]
The above isomorphism is given by a product decomposition
\[
X^\circ=X(1)\overset{D}{\rightsquigarrow} X\sqcup S^3(2).
\]

Under the product decomposition, choose a closed representative $G_X\subset\operatorname{int}X$
of $[G^{(m)}]$, and let $G_m^\circ$
be its image in $X^\circ \setminus\!\setminus D$ after reversing the product decomposition.

Now, let's check the Chern evaluation. Let $\mathbf y\in SFC(X^\circ,\gamma^\circ)$
corresponds to $\mathbf x\otimes\theta\in SFC(X,\mu_{K'})\otimes SFC(S^3(2))$ under the puncturing
identification. Using the adapted diagram of Lemma 9.13 in [3],
the flow lines selected by $\mathbf y$ avoid $D$; after decomposing along $D$,
those in the $X$-component are precisely the flow lines
selected by $\mathbf x$. Thus, by Definition 4.5 of [3], the vector field defining
$\mathfrak s(\mathbf y)$ restricts to a representative of
$\mathfrak s(\mathbf x)$ after an adjustment confined to
a collar of the new boundary disk (here, the boundary disk means the copy of $D$ exposed
on the $X$-component by the decomposition).
Choose this collar disjoint from $G_X$. Then, at each point of $G_m^\circ$, which is simply an unchanged copy of $G_X,$ the vector fields \(v_{\mathbf x}\) and \(v_{\mathbf y}\) give the same vector. Using the same metric, the planes perpendicular to these vectors are therefore also the same. Consequently,
\[
\left\langle c_1(\mathfrak s(\mathbf y)),[G_m^\circ]\right\rangle
=
\left\langle c_1(\mathfrak s(\mathbf x)),[G_X]\right\rangle,
\]
independently of $\theta$.
Since $G_m^\circ$ and $G^{(m)}$ become homologous after
filling the puncture, their classes in $H_2(X^\circ;\mathbb Z)$
differ by a multiple of the puncture sphere $S_0$.
The positive and negative regions of $S_0$ are both disks,
so $\langle c_1(\mathfrak s(\mathbf y)),[S_0]\rangle=1-1=0$ (See Proposition 3.4 of [4]). Hence
\[
\langle c_1(\mathfrak s(\mathbf y)),[G_m^\circ]\rangle
=
\langle c_1(\mathfrak s(\mathbf y)),[G^{(m)}]\rangle.\]
Therefore
\[
SFH(X^\circ,\gamma^\circ\mid G^{(m)})
\cong
SFH(X,\mu_{K'}\mid G^{(m)})\otimes V.
\]
Combining this with the knot-exterior identification yields
\[
SFH(X^\circ,\gamma^\circ\mid G^{(m)})
\cong
\widehat{HFK}
(Y',K'\mid G^{(m)})\otimes V.
\]

Hence, we have a clear relationship between $SFH(X^\circ,\gamma^\circ\mid G^{(m)})$ and $
\widehat{HFK}
(Y',K'\mid G^{(m)})$. It remains to compare $SFH(X^\circ,\gamma^\circ\mid G^{(m)})$ with $SFH(M_0(K'),\gamma_0(K'))$. A useful observation is that, for $\mathfrak{s}\in \rm Spin^c$$(X^\circ,\gamma^\circ |G^{(m)})$, we obtain a nice
splitting of the generators of the Heegaard Floer complex into intersections in $\widetilde{F}^\circ =\widetilde{F}\setminus (D_{z_k}\cup D_{w_k})\subset \widetilde{\Sigma}_{K'}^\circ$ and $\overline{F}^\circ=\overline{F}\setminus D_{w_r}.$

Suppose $\boldsymbol{\alpha}_{K'}=\boldsymbol{\alpha}^{K'}_{\text{left}}\cup\boldsymbol{\alpha}_{\text{right}},$ $\boldsymbol{\beta}_{K'}=\boldsymbol{\beta}^{K'}_{\text{left}}\cup\boldsymbol{\beta}_{\text{right}}.$ Let $\mathbb{T}_{{\alpha}_{K'}}=\prod_{\alpha\in \boldsymbol{{\alpha}}_{K'}}\alpha$ and $\mathbb{T}_{{\beta}_{K'}}=\prod_{\beta\in \boldsymbol{\beta}_{K'}}\beta.$ Let $C_{\text{left}}$ be the $\mathbb{F}_2$ vector space generated by the intersection points of $\mathbb{T}_{\alpha_{\text{left}}^{K'}}=\prod_{\alpha\in \boldsymbol{\alpha}^{K'}_{\text{left}}}\alpha$ and $\mathbb{T}_{{\beta}^{K'}_{\text{left}}}=\prod_{\beta\in \boldsymbol\beta^{K'}_\text{left}}\beta$ in $\operatorname{Sym}^{d'+1}(\widetilde{\Sigma}^\circ_{K'})$. Let $C_{\text{right}}$ be the $\mathbb{F}_2$ vector space generated by the intersection points of $\mathbb{T}_{\alpha_{\text{right}}}=\prod_{\alpha\in \boldsymbol{{\alpha}_{\text{right}}}}\alpha$ and $\mathbb{T}_{{\beta}_{\text{right}}}=\prod_{\beta\in \boldsymbol{{\beta}_\text{right}}}\beta$ in $\operatorname{Sym}^{2g+1}(\widetilde{\Sigma}^\circ_{K'})$.

\begin{lemma}
An intersection point $\mathbf{x}\in \mathbb{T}_{\alpha_{K'}}\cap\mathbb{T}_{\beta_{K'}}$ induces a $\text{Spin}^c$ structure $\mathfrak{s}(\mathbf{x})$ in $\rm Spin^c$$ (X^\circ,\gamma^\circ |G^{(m)})$ if and only if $\mathbf{x}\in C_{\rm{left}}\otimes C_{\rm{right}}$.
\end{lemma}
\begin{proof}
    By the index formula from [4], [5] and [6], $$\langle c_1(\mathfrak{s}(\mathbf{x})),H(\widetilde{F}^\circ)\rangle=e(\widetilde{F}^\circ)+2n_{\mathbf{x}}(\widetilde{F}^\circ).$$

    Since $e(\widetilde{F}^\circ)=-2g-2d'-2$, $H(\widetilde{F}^\circ)=[G^{(m)}]$, and $\mathfrak{s}(\mathbf{x})\in \text{Spin}^c (X^\circ,\gamma^\circ |G^{(m)})$ if and only if $$\langle c_1(\mathfrak{s}(\mathbf{x})),H(\widetilde{F}^\circ)\rangle=\langle c_1(\mathfrak{s}(\mathbf{x})),[G^{(m)}]\rangle=2-2g,$$ we know that $\mathfrak{s}(\mathbf{x})\in\text{Spin}^c (X^\circ,\gamma^\circ |G^{(m)})$ if and only if $n_{\mathbf{x}}(\widetilde{F}^\circ)=d'+2.$

    There are $d'+1$ alpha-curves in $\widetilde{\mathcal{H}}^\circ_{K'}$ lying entirely in $\widetilde{F}^\circ$, which are precisely the curves in $\boldsymbol{\alpha}_{\text{left}}^{K'}$, and $\alpha_0,\beta_0$ are on the boundary $\widetilde{F}^\circ.$ Therefore, $n_{\mathbf{x}}(\widetilde{F}^\circ)=d'+2$ if and only if $\mathbf{x}\in C_{\text{left}}\otimes C_{\text{right}}$.
\end{proof}

Note that there is an obvious one-to-one correspondence between generators of our $SFC(M_0(K'),\gamma_0(K'))$ and generators of $C_{\text{left}}$. Hence, to compare $SFH(X^\circ,\gamma^\circ\mid G^{(m)})$ with $SFH(M_0(K'),\gamma_0(K'))$, we will prove that $SFH(X^\circ,\gamma^\circ\mid G^{(m)})$ is isomorphic to the Floer homology of the chain complex $(C_{\text{left}},\partial_{\text{left}})\otimes (V,\operatorname{id})$. Before stating the theorem, we shall first give a rigorous definition of $(C_{\text{left}},\partial_{\text{left}})$.

Similar to Lekili's definition of $QFH'$ in [7], We may define $QKFH'(Y',K',f)$. Let $\widetilde{\Sigma}_{\max}^\circ=\widetilde{\Sigma}_{\max}\setminus (D_{z_k}\cup D_{w_k})$. Define the chain complex to be the $\mathbb{F}_2$ vector space generated by the intersection points of $\mathbb{T}_{\alpha_{\text{left}}^{K'}}$ and $\mathbb{T}_{{\beta}^{K'}_{\text{left}}}$ in $\operatorname{Sym}^{d'+1}(\widetilde{\Sigma}_{\max}^\circ)$, equipped with the diffential:

$$\partial_{K'}\mathbf{x}=\sum_{\varphi\in \pi^\circ_2(\mathbf{x},\mathbf{y}),\mu(\varphi)=1}\#\widehat{\mathcal{M}}(\varphi)\mathbf{y}.$$

There are two issues we need to address to make $QKFH'(Y',K',f)$ well-defined. The first issue is the compactness of the moduli space $\mathcal{M}(\varphi).$ The second issue is proving that $\partial_{K'}^2=0.$ To ensure compactness of the moduli space, we will adopt the cylindrical setting of Lipshitz in [5], where one considers pseudo-holomorphic curves in $\widetilde{\Sigma}_{\max}^\circ\times [0,1]\times \mathbb{R}$ instead of disks in $\operatorname{Sym}^{d'+1}(\widetilde{\Sigma}_{\max}^\circ)$, and choose an almost-complex structure $J$ on $\widetilde{\Sigma}_{\max}^\circ\times [0,1]\times \mathbb{R}$ that satisfies Lipshitz's conditions $\mathbf{J1}$ to $\mathbf{J5}$ [5]. Conditions
\(\mathbf{J3}-\mathbf{J5}\) imply that
$J=j_s\oplus j_D$
for a path \(j_s\), \(s\in[0,1]\), of complex structures on
\(\widetilde\Sigma_{\max}^{\circ}\). In the cylindrical setting, \(\pi_2^\circ(\mathbf x,\mathbf y)\) denote the set of relative
homology classes represented by proper maps
$
u:S\longrightarrow
\widetilde\Sigma_{\max}^{\circ}\times[0,1]\times\mathbb R
$
from a punctured bordered surface, with alpha-boundary on
$
\boldsymbol\alpha_{\mathrm{left}}^{K'}
\times\{1\}\times\mathbb R,
$
beta-boundary on
$
\boldsymbol\beta_{\mathrm{left}}^{K'}
\times\{0\}\times\mathbb R,
$
and asymptotic to \(\mathbf x\) and \(\mathbf y\) as
\(t\to-\infty\) and \(t\to+\infty\), respectively. We work in this cylindrical formulation
throughout the remainder of the section.

If \(v\) is an interior sphere bubble, then its projected domain \(\mathcal D(v)\) is a two-cycle on \(\widetilde {\Sigma}_{\max}^\circ\).  Since its
coefficient on regions adjacent to the boundary is zero and
\(\widetilde{\Sigma}_{\max}^\circ\) is connected, it follows that
$
\mathcal D(v)=0.$
The energy estimate therefore gives \(E(v)=0\), so \(v\) is
constant. Hence, bubbling at interior points can be ruled out.

The same zero-multiplicity condition also prevents escape toward
\(\partial\widetilde{\Sigma}_{\max}^\circ\).  Indeed, for \(p\in U_{\partial}
\cong
\partial\widetilde\Sigma_{\max}^{\circ}\times[0,\varepsilon)\), the vertical
cylinder
\[
Z_p=\{p\}\times[0,1]\times\mathbb R
\]
is \(J\)-holomorphic, and every counted curve \(u\) satisfies
\[
u\cdot Z_p=n_p(u)=0.
\]
No component of \(u\) is contained in \(Z_p\), and every isolated
intersection of two \(J\)-holomorphic curves has positive local
multiplicity.  Hence \(u\cap Z_p=\varnothing\).  Varying \(p\) over
\(U_{\partial}\) gives
\[
\pi_{\widetilde\Sigma_{\max}^{\circ}}(u)\cap U_{\partial}
=\varnothing,
\]
so the curves cannot escape through the puncture ends.

Furthermore, the attaching curves $\boldsymbol{\alpha}_{\text{left}}^s$ are linearly independent in
\(H_1(\widetilde{\Sigma}_{\max}^\circ)\) and so are the curves $\boldsymbol{\beta}_{\text{left}}^s$. The only additional
relations in the enlarged systems are generated by
\(\alpha_r\), respectively \(\beta_r\).

\begin{figure}[H]
 \centering
  \includegraphics[width=0.5\textwidth]{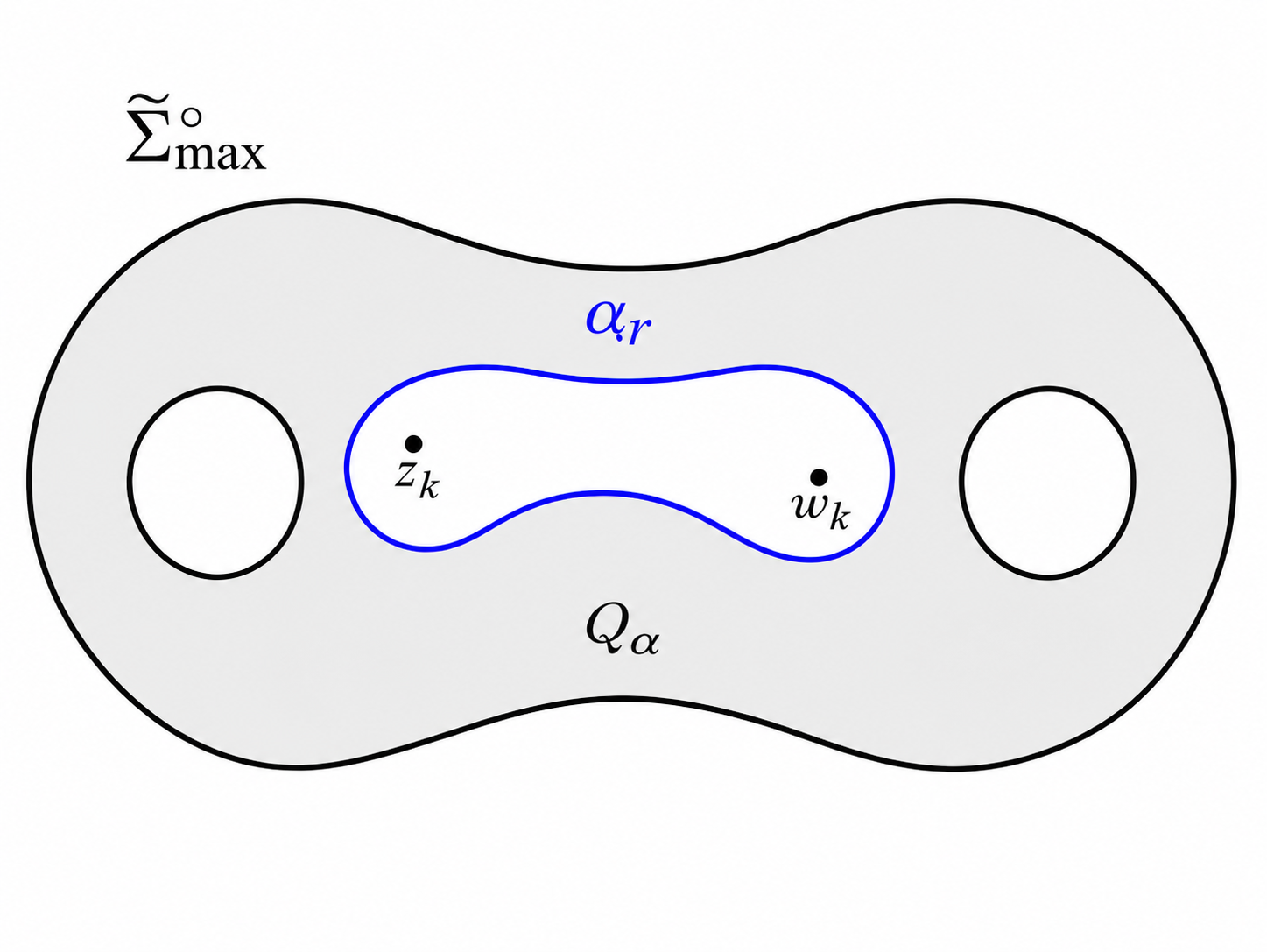}
 \caption{The region $Q_\alpha$.}
 \end{figure}
As shown in Figure 4, let \(Q_\alpha\) be the component of
\(\widetilde{\Sigma}_{\max}^\circ\setminus\alpha_r\) that is disjoint from
\(\partial\widetilde{\Sigma}_{\max}^\circ\), oriented so that
\[
\partial Q_\alpha=\alpha_r.
\]
Define \(Q_\beta\) analogously.  Notice that \(\alpha_r\) and
\(\beta_r\) are null-homologous in \(\widetilde{\Sigma}_{\max}^\circ\), although
they are not null-homotopic: the other side of each curve is a pair
of pants containing the two boundary components arising from
\(w_K\) and \(z_K\).

Suppose \(A\) is the domain of an alpha-boundary degeneration, its projected domain satisfies
\[
\partial A
=
\sum_{\alpha_t\in\boldsymbol{\alpha}_{\text{left}}^s}
m_{\alpha_t}\alpha_t
+
m_{\alpha_r}\alpha_r
\]
for some integers \(a_i,m\).  Since \(A\) is a two-chain, and $\alpha_r$ is null-homologous,
\[
0=[\partial A]
=
\sum_{\alpha_t\in\boldsymbol{\alpha}_{\text{left}}^s}
m_{\alpha_t}[\alpha_t]
+
m_{\alpha_r}[\alpha_r]
=
\sum_{\alpha_t\in\boldsymbol{\alpha}_{\text{left}}^s}
m_{\alpha_t}[\alpha_t].
\]
The classes represented by the curves in
\(\boldsymbol{\alpha}_{\text{left}}^s\) are linearly independent, and hence each
\(m_{\alpha_t}=0\).  It follows that
\[
\partial A=m_{\alpha_r}\alpha_r.
\]

Consequently,
\[
\partial(A-m_{\alpha_r}Q_\alpha)=0.
\]
A domain on \(\widetilde{\Sigma}_{\max}^\circ\) has multiplicity zero in every
region meeting \(\partial\widetilde{\Sigma}_{\max}^\circ\).  On the other hand, a
domain with zero formal boundary has equal multiplicities in any
two adjacent regions.  Since the dual graph of the regions is
connected, all coefficients of \(A-m_{\alpha_r}Q_\alpha\) are therefore
equal; the boundary normalization forces this common coefficient
to be zero.  Hence
\[
A=m_{\alpha_r}Q_\alpha.
\]
Positivity of domains gives \(m_{\alpha_r}\geq0\), with \(m_{\alpha_r}>0\) whenever the
degeneration is nonconstant. The beta-boundary case is identical.

Choose disjoint separating curves
\[
c_1^\alpha,\ldots,c_{d'+1}^\alpha
\subset\widetilde{\Sigma}_{\max}^\circ
\]
which are disjoint from all the alpha-curves and such that cutting
along them decomposes \(\widetilde{\Sigma}_{\max}^\circ\) into
\[
H_1^\alpha,\ldots,H_{d'}^\alpha,R_\alpha,
\]
where each \(H_i^\alpha\) is a one-holed torus containing
a distinct curve in $\boldsymbol{\alpha}_{\text{left}}^s$, while \(R_\alpha\) has genus \(g\), contains
\(\alpha_r\), and contains the two components of
\(\partial\widetilde{\Sigma}_{\max}^\circ\).

The curve \(\alpha_r\) is homotopically essential in \(R_\alpha\).

Cut \(\widetilde{\Sigma}_{\max}^\circ\) along the curves \(c_i^\alpha\).  Fix
orientation-compatible complex structures on the resulting
components
\[
H_1^\alpha,\ldots,H_{d'}^\alpha,R_\alpha
\]
which are cylindrical near the newly created boundary circles. The complex structures are also fixed near the original boundary
\(\partial\widetilde{\Sigma}_{\max}^\circ\) and near all the alpha-curves.

Cutting along each $c_i^\alpha$, creates a pair of new boundary circles. For \(T>0\), glue each such pair of new boundary circles
by the standard complex cylinder. The resulting Riemann surface is diffeomorphic to
\(\widetilde{\Sigma}_{\max}^\circ\).  Pull the resulting complex structure back to $\widetilde{{\Sigma}}^\circ_{\max}$ via a diffeomorphism which is fixed away
from prescribed annular neighborhoods of the \(c_i^\alpha\), and
denote this complex structure on $\widetilde{{\Sigma}}^\circ_{\max}$ by \(j_{\alpha,T}\).

Thus \(j_{\alpha,T}\) is independent of \(T\) near
\(\partial\widetilde{\Sigma}_{\max}^\circ\) and near the alpha-curves, and the
annular neighborhood of each \(c_i^\alpha\) is a conformal
cylinder of length \(2T\).  As \(T\to\infty\), these cylinders
break into pairs of semi-infinite cylindrical ends; equivalently,
after compactifying the ends, $(\widetilde{\Sigma}_{\max}^\circ, j_{\alpha,T})$ converges to the nodal
union
\[
H_1^\alpha\cup\cdots\cup H_{d'}^\alpha\cup R_\alpha.
\]

At the
alpha-boundary of the cylindrical target $\widetilde{\Sigma}_{\max}^\circ\times [0,1]\times \mathbb{R}$, pick the complex structure $J_T$ on the target such that
\[
J_T\big|_{\widetilde{\Sigma}_{\max}^\circ\times\{1\}\times\mathbb R}
=
j_{\alpha,T}\oplus j_D,
\]
where \(j_D\) is the standard complex structure on the strip $[0,1]\times\mathbb{R}$.

Fix \(m>0\).  We claim that there is \(T_0(m)>0\) such that, for
every \(T\geq T_0(m)\), there is no nonconstant
\(J_T\)-holomorphic alpha-boundary degeneration with domain
\(mQ_\alpha\).

Suppose otherwise. Then is a sequence $T_n\longrightarrow\infty$ and, for every \(n\), there is a \(J_{T_n}\)-holomorphic alpha-boundary
degeneration
$
u_n$ with domain $mQ_\alpha.$

Because \(J_{T_n}=j_{\alpha,T_n}\oplus j_D\) near the
alpha-boundary, the cylindrical boundary degeneration \(u_n\)
determines a
\(\operatorname{Sym}^{d'+1}(j_{\alpha,T_n})\)-holomorphic disk
\[
v_n:
(D,\partial D)
\longrightarrow
\left(
\operatorname{Sym}^{d'+1}(\widetilde{\Sigma}_{\max}^\circ),
\mathbb T_{\alpha_{\text{left}}^{K'}}
\right).
\]
Indeed, the projection of \(u_n\) to
the strip is an $(d'+1)$-fold branched covering. For a point \(q\) in the strip, collecting
the images in \(\widetilde{\Sigma}_{\max}^\circ\) of the \(d'+1\) points lying over
\(q\) gives an unordered divisor
\[
v_n(q)\in\operatorname{Sym}^{d'+1}(\widetilde{\Sigma}_{\max}^\circ).
\]
Because \(J_T=j_{\alpha,T_n}\oplus j_D\) at this boundary, the
surface and strip components of the Cauchy--Riemann equation
separate, and the divisor-valued map \(v_n\) is holomorphic.

After passing to a subsequence, neck-stretching compactness gives a
stable holomorphic limit $v_\infty$
in the symmetric product of the nodal surface obtained as
\(T_n\to\infty\).  The stable map \(v_\infty\) may a priori have
several components, and its irreducible component containing the limiting
alpha-torus consists of divisors having one point in each of $H_1^\alpha,\cdots H_{d'}^\alpha,R_\alpha$, and is consequently
\[
X_\alpha
=
H_1^\alpha
\times\cdots\times
H_{d'}^\alpha
\times R_\alpha.
\]
The limiting boundary condition is
\[
L_\alpha
=
\alpha_1
\times\cdots\times
\alpha_{d'}
\times\alpha_r,
\]
where each $\alpha_i$ is the distinct curve in $\boldsymbol{\alpha}_{\text{left}}^s$ corresponding to $H_i^\alpha.$

For each \(i\), the surface \(H_i^\alpha\) has positive genus and
\(\alpha_i\) is essential in \(H_i^\alpha\).  The homotopy exact
sequence of the pair gives
\[
\pi_2(H_i^\alpha,\alpha_i)=0.
\]
Likewise,
\[
\pi_2(R_\alpha,\alpha_r)=0,
\]
since \(R_\alpha\) is aspherical and \(\alpha_r\) is essential in
\(R_\alpha\).

Maps of pairs into a product split coordinatewise. Hence
\[
\pi_2(X_\alpha,L_\alpha)
\cong
\prod_{i=1}^{d'}
\pi_2(H_i^\alpha,\alpha_i)
\times
\pi_2(R_\alpha,\alpha_r)
=
0.
\]
The boundary-bearing disk component of the limiting stable map is
therefore relatively null-homotopic. Its symplectic area is zero,
so it is constant.

The constant disk maps to the limiting generator, which is
disjoint from the degeneration nodes.  Moreover,
\(\pi_2(X_\alpha)=0\), so no sphere can bubble inside
\(X_\alpha\); components lying in other irreducible components of
the nodal symmetric product can attach only along the nodal loci.
It follows that the entire stable limit is constant.

This contradicts the assumption that the original domains are
\(mQ_\alpha\). Indeed, if
\[
p\in\operatorname{int}Q_\alpha
\]
is chosen away from the stretching curves, then
\[
n_p(u_n)=m>0,
\]
whereas the constant limiting map has multiplicity zero at \(p\).
Thus, for sufficiently large \(T\), there is no alpha-boundary
degeneration with domain \(mQ_\alpha\).

By admissibility, for each pair of generators there are only
finitely many nonnegative connecting domains.  Since the set of
generators is finite, only finitely many index-one and index-two
total domains occur. For each such domain \(A\), if a stable degeneration in class \(A\)
contains an alpha-boundary component with domain \(mQ_\alpha\),
then, for any $p\in\operatorname{int}Q_\alpha$
away from the stretching regions, we have $n_p(A)
\geq m$, and hence \(m\) is bounded. Thus only finitely many
multiplicities \(m\) can occur, so one sufficiently large neck which excludes
all alpha-boundary degenerations relevant to the moduli spaces
under consideration.

The beta-boundary case is identical. Now, choose \(T_\alpha,T_\beta\) sufficiently large so that the preceding
alpha- and beta-boundary-degeneration exclusions hold.  Since
\(j_{\alpha,T_\alpha}\) and \(j_{\beta,T_\beta}\) agree with the
same fixed models near the puncture ends and near the points used
in \(\mathbf{J2}\), choose a generic path \(j_s\), fixed on these
neighborhoods and constant near its endpoints, such that
\[
j_s=j_{\beta,T_\beta}
\quad\text{for \(s\) near \(0\),}
\qquad
j_s=j_{\alpha,T_\alpha}
\quad\text{for \(s\) near \(1\).}
\]
Set
\[
J=j_s\oplus j_D.
\]

Hence, we can also avoid boundary bubbles by choosing an appropriate almost complex structure, and the compactness of $\mathcal{M}(\varphi)$ is ensured.

Since \(J\) satisfies \(\mathbf{J5}\), the strip projection
\[
\pi_D:
\widetilde\Sigma_{\max}^{\circ}\times [0,1]\times\mathbb{R}
\longrightarrow [0,1]\times\mathbb{R}
\]
is \(J\)-holomorphic. Hence Lipshitz's tautological
correspondence [5] identifies the cylindrical moduli spaces with the
corresponding Whitney-disk moduli spaces, calculated with respect to the induced path
$
\operatorname{Sym}^{d'+1}(j_s)$
of complex structures on
\(\operatorname{Sym}^{d'+1}
(\widetilde\Sigma_{\max}^{\circ})\).
Indeed, every relevant domain has coefficient zero in the regions
meeting
\(\partial\widetilde\Sigma_{\max}^{\circ}\), and the preceding
vertical-cylinder argument shows that every counted curve remains
in a compact subset of the interior. Thus the cylindrical and
Whitney-disk counts agree.

The verification of $\partial_{K'}^2=0$ will be postponed and will be part of the proof of the identification between $QKFH'(Y',K',f)\otimes V$ and $SFH(X^\circ,\gamma^\circ|G^{(m)})$.

Now we turn to the comparison between $QKFH'(Y',K',f)\otimes V$ and $SFH(X^\circ,\gamma^\circ|G^{(m)})$. We need the following two lemmas:
\begin{lemma}
\label{lem:knot-Lekili-12}
Let \(\alpha'_0\) and \(\beta'_0\) be parallel pushoffs of \(\alpha_0\) and
\(\beta_0\) into the interior \(\widetilde F^\circ\). Suppose the attaching curves divide the cylindrical region between \(\alpha'_0\) and \(\alpha_0\) into components $
a_1,\cdots,a_{4g}$ and divide the region between $\beta'_0$ and \(\beta_0\) into $b_1,\cdots, b_{4g}$.
Choose the labeling so that $a_1$ and $b_1$ lie in the region containing $\partial D_{w_k}.$

Let $\mathbf x=\mathbf x_{\mathrm{left}}\otimes
\mathbf x_{\mathrm{right}}$ and $\mathbf y=\mathbf y_{\mathrm{left}}\otimes
\mathbf y_{\mathrm{right}}$ be in $C_{\mathrm{left}}\otimes C_{\mathrm{right}}$ and $A\in \pi_2(x,y)$ and let $u$ be a Maslov index 1 holomorphic curve in the homology class $A.$ Assume moreover that the contribution of curves in the class $A$ to the differential is nonzero. Then,
\[n_{a_1}(u)=\cdots=n_{a_{4g}}(u)=n_{b_1}(u)=\cdots=n_{b_{4g}}(u)=0.\]
Furthermore, if the projection of the image of $u$ to the Heegaard surface $\widetilde{\Sigma}_{K'}^\circ =\widetilde{F}^\circ\cup_{\alpha_0\cup \beta_0}\overline{F}^\circ$ lies entirely in $\widetilde{F}^\circ$, one can find almost complex structures $j_0$ and $j_1$ on $\widetilde{\Sigma}_*$ and an admissible almost complex structure $J$ on $\widetilde{\Sigma}_{K'}^\circ \times [0,1]\times \mathbb{R}$ such that $J|_{\widetilde{\Sigma}_{K'}^\circ}\times \{0\}\times \mathbb{R}=j_0$ and $J|_{\widetilde{\Sigma}_{K'}^\circ}\times \{1\}\times \mathbb{R}=j_1$ with the property that $u$ restricted to the boundary does not hit $\{\alpha'_0,\beta'_0\}\times \{0,1\}\times \mathbb{R}.$
\end{lemma}

\begin{proof}
Suppose that
\[
n_{a_i}(u)\neq n_{a_{i+1}}(u)
\]
for some \(i\), and let \(\beta_j\) be the beta-arc separating
\(A_i\) from \(A_{i+1}\).  Stretching along \(\alpha'_0\) and \(\beta'_0\). For a neck-stretching family of complex structures \(\{J_T\}\), choose \(T_n\to\infty\) and suppose for each $n$ we have a \(J_{T_n}\)-holomorphic representatives \(u_n\) of \(A\). Lekili's argument in Lemma 12 of [7] applies verbatim up to the following conclusion:
after passing to a convergent subsequence, the neck-stretching limit of the curves the curves \(u_n\) is a pair of height-one holomorphic combs \(u_0\) and \(u_1\), and the multiplicity jump forces \(u_0\) to have a boundary puncture asymptotic to a Reeb chord on \(\alpha'_0\); let \(S_0\) be the connected component of the source of \(u_0\) containing this puncture, and then \(u_0(S_0)\) defines a two-chain $P$ in \(\widetilde{\Sigma}_{\max}^\circ\) whose boundary is supported
entirely on beta-curves.

Let \(\widehat\beta_j\) denote the closed curve obtained from
\(\beta_j\) after \(\alpha'_0\) is pinched.  We may write
\[
\partial P
=
m_{\widehat{\beta}_j}\widehat\beta_j
+\sum_{\beta\in\mathcal B_j}m_\beta\beta
+m_{\beta_r}\beta_r,
\]
where \(\mathcal B_j\) consists of the remaining
beta-curves aside from $\beta_r,$ $m_{\widehat{\beta}_j}=
\pm\bigl(n_{A_i}(u)-n_{A_{i+1}}(u)\bigr)\neq0.$
Here \(\widehat\beta_j\neq\beta_r\), since \(K'\subset\operatorname{int}(M_0)\), which means its bridge
presentation and the disk \(E_\beta^r\) defining
\(\beta_r\) may be chosen in the interior of the
maximal-fiber region, away from the seam \(\alpha_0\), and thus $\alpha_0\cap\beta_r=\varnothing.$

Since \(P\) is a two-chain, $[\partial P]=0.$
Moreover, \([\beta_r]=0\) because \(\beta_r\) bounds
\(Q_\beta\).  Hence
\[
m_{\widehat{\beta}_j}[\widehat\beta_j]
+\sum_{\beta\in\mathcal B_j}m_\beta [\beta]
=0.
\]
The non-\(r\) beta-curves are linearly independent, so
\(m_{\widehat{\beta}_j}=0\), a contradiction.  Therefore all \(a_i\)-multiplicities
are equal. The same argument with alpha and beta interchanged
proves the corresponding equality for the \(b_i\).

Finally, \(a_1\) and \(b_1\) lie in a region meeting
\(\partial D_{w_K}\), whose multiplicity is zero.  Thus all the
displayed multiplicities vanish.

For the second assertion, suppose that the boundary of a
left-supported curve continues to meet \(\alpha'_0\) for arbitrarily long
seam necks.  Following Lekili, we stretch along \(\alpha'_0\) and \(\beta'_0\),
pass to the limiting holomorphic comb, and rescale the component
whose projection to the strip is constant.  After a further
degeneration of the maximal-fiber side, this produces an
\((d'+2)\)-fold branched covering to the upper half-plane
\[
\coprod_V S_\nu\longrightarrow\mathbb H.
\]

Only the auxiliary degeneration requires modification. Choose $\{ S_v\}$
so that each left beta-curve aside from $\beta_r$ lies in a separate one-holed
tori in $\{S_v\}$, \(\beta_r\) lies in a planar piece in $\{S_v\}$ containing the two knot
punctures, and the closed seam curve \(\widehat\beta_j\) lies in the
remaining positive-genus piece. Each component of a covering of the connected space \(\mathbb H\) has positive integer degree and their degrees sum to \(d'+2\), each of the \(d'+2\) components would have degree \(1\). Then, some positive genus surface in $\{S_v\}$ would be
biholomorphic to \(\mathbb H\), which is a contradiction since $\mathbb{H}$ has genus 0. Similarly, the boundary of a
left-supported curve cannot meet $\beta'_0.$
\end{proof}
\begin{figure}[H]
 \centering
  \includegraphics[width=0.6\textwidth]{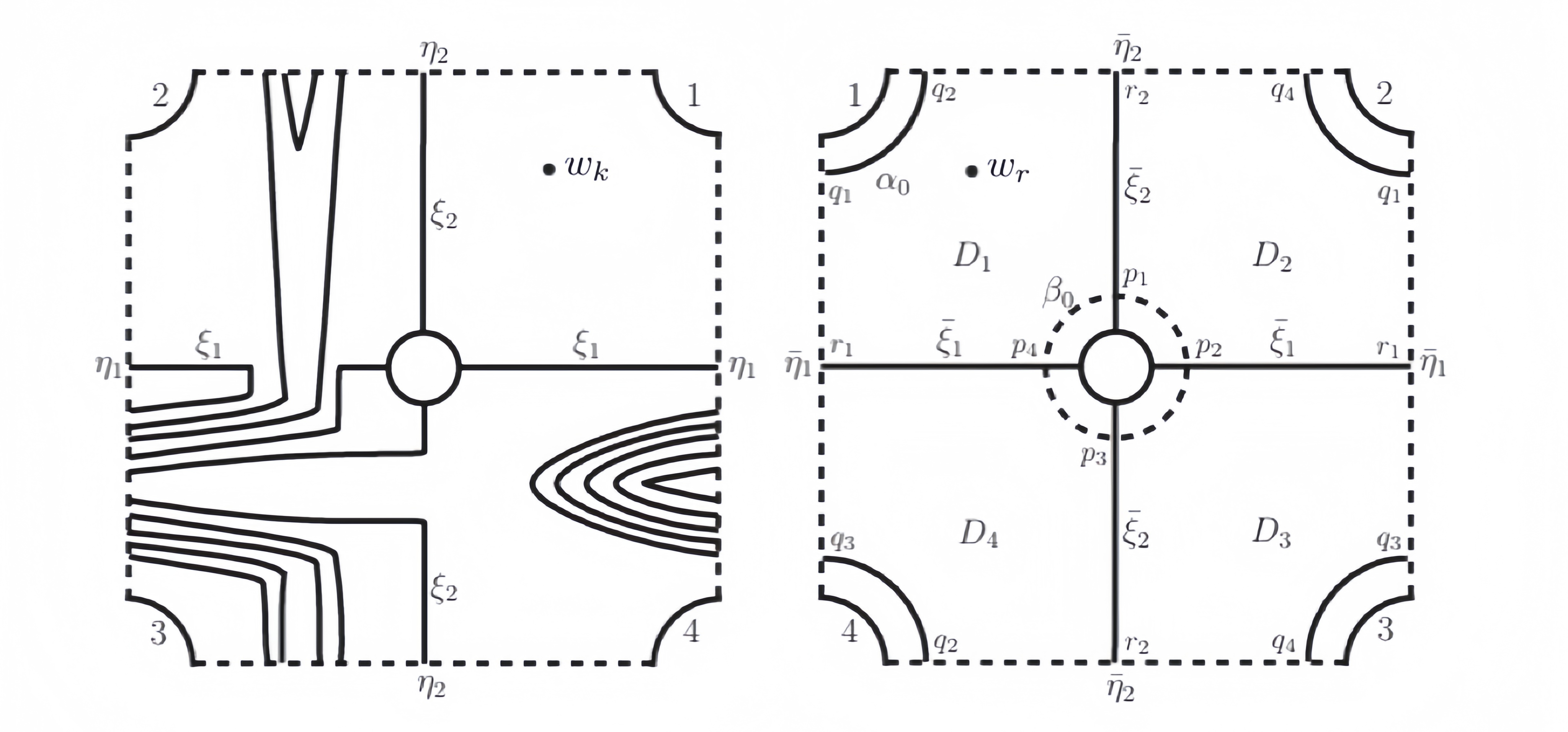}
 \caption{Example diagam for torus bundles.}
 \end{figure}
\begin{lemma}
As shown in Figure 5, the right side of the closed diagram $\overline{F}$ consists of hexagon regions $D_1,\ldots,D_{4g}$. For each \(i\), we let
\(\mathbf x_i,\mathbf y_i\in C_{\mathrm{ right}}\) be the pair of
right-hand generators for which $D_i$ induces a Whitney disk in $\operatorname{Sym}^{d+2g+1}(\Sigma_*)$
connecting $\mathbf x_L\otimes\mathbf x_i$ to
$\mathbf x_L\otimes\mathbf y_i$, with boundary on
$\mathbb T_{\alpha_*}=\prod_{\alpha\in\alpha_*}\alpha$ and $\mathbb T_{\beta_*}=\prod_{\beta\in\beta_*}\beta$. Here, $\mathbf x_L$ is a fixed left-hand generator which occupies each curve in $\boldsymbol{\alpha}_{\mathrm{left}}$ and
each curve in $\boldsymbol{\beta}_{\mathrm{left}}$ exactly once. In particular, we suppose $D_1 $ contains $\partial D_{w_r}$.

Consider $\mathbf{x}_{\mathrm{left}}\otimes \mathbf{x}_{\mathrm{right}},\mathbf{y}_{\mathrm{left}}\otimes \mathbf{y}_{\mathrm{right}}\in SFC(X^\circ,\gamma^\circ)$, with $\mathbf{x}_{\mathrm{right}}=\mathbf{x}_1$ or $\mathbf{y}_1$. Let $\varphi\in \pi_2(\mathbf{x}_{\mathrm{left}}\otimes \mathbf{x}_{\mathrm{right}},\mathbf{y}_{\mathrm{left}}\otimes \mathbf{y}_{\mathrm{right}})$ be a class that contributes non-trivially to the differential for given $\mathbf x_{\mathrm{left}},\mathbf{x}_{\mathrm{right}},\mathbf y_{\mathrm{left}},\mathbf y_{\mathrm{right}}$ and $u$ be a cylindrical $J$-holomorphic representative of
$\varphi$.
Then $\mathbf y_{\mathrm{right}}= \mathbf{x}_{\mathrm{right}},$
and the domain of $u$ is contained on the left side of the Heegaard diagram.
\end{lemma}

\begin{proof}
The case when $\mathbf{x}_{\mathrm{right}}=\mathbf{x}_1$ is just the same proof as Lemma 15 of [7].

Now, suppose $\mathbf{x}_{\mathrm{right}}=\mathbf{y}_1$.

By Lemma 4.4, the domain of \(u\) has
multiplicity zero in the annuli between
\(\alpha'_0\) and \(\alpha_0\), and between \(\beta'_0\) and \(\beta_0\).
Indeed, for each point \(p\) in these annuli,
\[
n_p(u)=u\cdot
\bigl(\{p\}\times[0,1]\times\mathbb R\bigr)=0.
\]
Since \(p\) is disjoint from the attaching
curves, \(u\) and \(Z_p\) cannot share a component. Hence every
intersection of \(u\) with \(Z_p\) is isolated and hence has positive local multiplicity;
in particular, such intersections cannot cancel algebraically.
Therefore
\[
u\cdot Z_p=n_p(u)=0
\]
implies \(u\cap Z_p=\varnothing\).  Varying \(p\) over the two
annuli shows that the surface projection of \(u\) avoids both
annuli. Hence stretching
along \(\alpha'_0\) and \(\beta'_0\) produces no nontrivial Reeb-chord ends at
the seams. Denote the left- and right-hand parts of the resulting
limiting building by $u_L$ and $u_R$, respectively.
The right-hand part $u_R$ is thus uncoupled from $u_L$
and has projected domain
\[
A_R:=\mathcal D(u_R)
=\mathcal D(u)|_{\overline F^\circ},
\]
which is a nonnegative domain from $\mathbf y_1$ to
the outgoing right-hand generator $\mathbf z_R$.

We first classify the
possible endpoints of \(A_{ R}\). Following discussion in Section 2.3 of [7], there is no right-supported class from
\(\mathbf x_i\) to \(\mathbf y_j\) when \(i\neq j\).  It follows,
by reversing \(A_{R}\), and, when necessary, composing with
the hexagonal class \(D_j\) from \(\mathbf x_j\) to
\(\mathbf y_j\), that
\[
\mathbf z_{ R}\in\{\mathbf y_1,\mathbf x_1\}.
\]

If \(\mathbf z_{\mathrm R}=\mathbf y_1\), then \(A_{\mathrm R}\)
has no corners on the right. Since its coefficient in the region
containing \(D_{w_r}\) is zero, it extends over \(D_{w_r}\) and,
after reflection across the two seams, determines a periodic domain
\(P\) in the standard diagram for \(S^1\times G^{(m)}\), with
$
P|_{\overline F}=A_{ R}.
$ By Lemma 4.4, the multiplicities of \(P\)
in all regions immediately to the left of \(\alpha_0\) and
\(\beta_0\) are equal. By observation about the periodic domains in the standard diagram for \(S^1\times G^{(m)}\), we see that any combination of periodic domains in \(S^1\times G^{(m)}\) with such property must have trivial restriction to $\overline{F}$, it follows that
\[
A_{ R}=P|_{\overline F}=0.
\]

Suppose instead that $\mathbf z_{ R}=\mathbf x_1.$ Extend $A_R$ over $D_{w_r}$ with coefficient zero,
and denote the resulting domain on $\overline F$
by $\overline A_R$.

Both $\overline A_R$ and $\overline F-D_1$ connect
$\mathbf y_1$ to $\mathbf x_1$ on the right side. Hence, their difference $Q=\overline A_R-(\overline F-D_1)$ has no right-hand corners and has coefficient zero
in $D_1$. Completing $Q$ by reflection across the
two seams gives a periodic domain $P$ in the standard
diagram for $S^1\times G^{(m)}$, with $P|_{\overline F}=Q.$
With the same argument as in the preceding
case, we have $P|_{\overline F}=0$. Therefore $\overline A_R=\overline F-D_1.$

Set $B_R=\overline F-D_1$. The domain $B_R$ is the extension of $A_R$ by zero over
$D_{w_r}$. Since it vanishes near this puncture, capping
does not change its Euler or point measures.
Lipshitz's cylindrical index formula therefore gives
\[
\operatorname{ind}(u_R)
=
e(B_R)+n_{\mathbf y_1}(B_R)+n_{\mathbf x_1}(B_R).
\]
Using
\[
e(\overline F)=-2g,\qquad
n_{\mathbf x_1}(\overline F)
=n_{\mathbf y_1}(\overline F)=2g,\qquad
\mu(D_1)=1,
\]
linearity yields
\[
\operatorname{ind}(u_R)
=
-2g+2g+2g-1
=
2g-1.
\]

There are no seam ends, so the indices of the uncoupled
left- and right-hand components add. Since the original
curve has index one, this would force
\[
\operatorname{ind}(u_L)
=
1-(2g-1)=2-2g<0,
\]
as $g\geq2$. This is impossible for a regular left-hand
component; trivial strips have index zero.
Thus $\mathbf z_R=\mathbf x_1$ is excluded.
The preceding case then gives $A_{ R}=0$ and $\mathbf z_{ R}=\mathbf y_1,$
as claimed.
\end{proof}
Using the above two lemmas, we can compute the ``right side homology" for $\widetilde{\mathcal{H}}_{K'}^\circ.$ That is, equip $C_{\mathrm{right}}$ with the differential $\partial_{\mathrm{right}}$ that counts Maslov-index-one holomorphic curves with domain supported entirely in the punctured right-hand region \(\overline{F}^\circ\); the homology group of the resulting chain complex $(C_{\mathrm{right}},\partial_{\mathrm{right}})$ is rank 2 with generators $\mathbf{x}_1,\mathbf{y}_1.$

We know that $SFC(X^\circ,\gamma^\circ|G^{(m)})=C_{\mathrm{left}}\otimes C_{\mathrm{right}}$. We denote by $\partial_{\mathrm{left}}\otimes\mathbbm{1}$ the contribution of those holomorphic curves whose domain lie in $\widetilde{F}^\circ$ and which act by the identify on $C_{\mathrm{right}}.$ Also, we denote by $\mathbbm{1}\otimes \partial_{\mathrm{right}}$ the contribution of those holomorphic curves whose domain lie in $\overline{F}^\circ$ and which act by the identify on $C_{\mathrm{left}}.$
\begin{proposition}
    $SFH(X^\circ,\gamma^\circ,\mathfrak{s})\simeq H(C_{\mathrm{left}}\otimes V,\partial_{\mathrm{left}}\otimes \mathbbm{1},\mathfrak{s})$ for $\mathfrak{s}\in \mathrm{Spin}^c(X^\circ,\gamma^\circ|G^{(m)}).$
\end{proposition}
\begin{proof}
    Let $V$ be generated by $\mathbf{x}_1,$ ${\mathbf{y}}_1$. We claim that there is a chain map $\iota:C_{\mathrm{left}}\otimes V\to C_{\mathrm{left}}\otimes C_{\mathrm{right}}$ given by natural inclusion. By Lemma 4.5, if a holomorphic map contributing to the differential originates at $[x_{\mathrm{left}}\otimes \mathbf{x}_{\mathrm{right}}]$, with $\mathbf{x}_1$ or $ \mathbf{y}_1$, then it has to converge to a generator of the form  $[\mathbf{y}_{\mathrm{left}}\otimes \mathbf{x}_{\mathrm{right}}]$, and the domain of the map has to lie on the left half of the Heegaard diagram; these are exactly the contributions to the differential captured by $\partial_{\mathrm{left}}\otimes \mathbbm {1}.$ Hence, $\iota$ is a chain map.

    Now we shall show that $\iota$ induces an isomorphism on the homology group.

    By Lemma 4.4, we can identify the differential on $C_{\mathrm{left}}\otimes C_{\mathrm{right}}\subset SFC(X^\circ,\gamma^\circ)$ with $\partial_{\mathrm{left}}\otimes\mathbbm{1}+\mathbbm{1}\otimes \partial_{\mathrm{right}}$. The injectivity of $\iota$ in homology follows from the fact that $\mathbf{x}_1,\mathbf{y}_1$ does not lie in the image of $\partial_{\mathrm{right}}$. Thus, we only need to check that $\iota$ is surjective in homology.

    Suppose that $\mathbf{s}_1\mathbf{x}_1+\cdots+\mathbf{s}_{4g}\mathbf{x}_{4g}+\mathbf{t}_1\mathbf{y}_1+\cdots+\mathbf{t}_{4g}\mathbf{y}_{4g}\in C_{\mathrm{left}}\otimes C_{\mathrm{right}}$ is in the kernel of $\partial_{\mathrm{left}}\otimes\mathbbm{1}+\mathbbm{1}\otimes \partial_{\mathrm{right}}$, where we have chosen the notation so that $\mathbf{s}_i,\mathbf{t}_i$ are elements in $C_{\mathrm{left}}$ and $\mathbf{x}_i,\mathbf{y}_i$ are generators in $C_{\mathrm{right}}$ we defined before. Following Lekili's computation in Section 2.3 of [7], we have $\partial_{\mathrm{right}}\mathbf{x}_1=\partial_{\mathrm{right}}\mathbf{y}_1=0$, and $\partial_{\mathrm{right}}\mathbf{x}_i=\mathbf{y}_i$ for $i\geq 2.$

    Now, because this element is in the kernel of $\partial_{\mathrm{left}}\otimes\mathbbm{1}+\mathbbm{1}\otimes \partial_{\mathrm{right}}$,
    \begin{align*}
        \partial_{\mathrm{left}}\mathbf{s}_1=0\qquad&\text{and}\qquad \partial_{\mathrm{left}}\mathbf{t}_1=0\\
        \partial_{\mathrm{left}}\mathbf{s}_i=0\qquad&\text{and}\qquad  \mathbf{s}_i+\partial_{\mathrm{left}}\mathbf{t}_i=0\quad \text{for $i>1$.}
    \end{align*}

Now, observe that the above equation gives $(\partial_{\mathrm{left}}\otimes\mathbbm{1}+\mathbbm{1}\otimes \partial_{\mathrm{right}})(\mathbf{t}_i\mathbf{x}_i)=-\mathbf{s}_i\mathbf{x}_i+\mathbf{t}_i\mathbf{y}_i$ for $i>1.$ Hence, $-\mathbf{s}_i\mathbf{x}_i+\mathbf{t}_i\mathbf{y}_i$ is in the kernel. Since we are working with $\mathbb{F}_2$ coefficients, $\mathbf{s}_i\mathbf{x}_i+\mathbf{t}_i\mathbf{y}_i=-\mathbf{s}_i\mathbf{x}_i+\mathbf{t}_i\mathbf{y}_i$ is in the kernel.

Therefore, we can ignore all the terms other than $\mathbf{s}_1\mathbf{x}_1+\mathbf{t}_1\mathbf{y}_1$, which is in the image of $\iota$ as desired.
\end{proof}
\begin{theorem}
There is an isomorphism $QKFH'(Y',K',f)\otimes V\cong SFH(X^\circ,\gamma^\circ\mid G^{(m)}).$
\end{theorem}

\begin{proof}
Let \(C_{K'}\) denote the \(\mathbb F_2\)-vector space used in the
above definition of \(QKFH'(Y',K',f)\). By Proposition 4.6, it suffices to identify
$(
C_{\mathrm{left}}\otimes \mathbbm{1},\,
\partial_{\mathrm{left}}\otimes\mathbbm{1})$
with
$(
C_{K'}\otimes V,\,
\partial_{K'}\otimes\mathbbm{1})$.

There is an evident one-to-one correspondence between their
generators. Choose the almost-complex structures in the two cylindrical targets
to agree on the common maximal-fiber region.  Fix
\(\theta\in\{\mathbf x_1,\mathbf y_1\}\).  By
Lemma 4.4, every curve contributing to
\(\partial_{\mathrm{left}}\otimes\mathbbm{1}\) has
multiplicity zero in the two seam annuli.  Positivity of
intersections therefore shows that it stays away from the seams.
Consequently, after removing the trivial right-hand strips over
\(\theta\), it is precisely a curve counted by
\(\partial_{K'}\).

Conversely, adjoining the trivial right-hand strips over \(\theta\)
to a curve counted by \(\partial_{K'}\) gives a curve counted by
\(\partial_{\mathrm{left}}\otimes\mathbbm{1}\).  These two
operations are inverse and preserve the Fredholm index.  Hence the
corresponding index-one moduli spaces have the same mod-two counts,
and therefore
$
\partial_{\mathrm{left}}=\partial_{K'}.$
The two displayed complexes are thus isomorphic.  Taking homology
and applying Proposition 4.6 proves the
result.
\end{proof}
\begin{theorem}
    $SFH(M_0(K'),\gamma_0(K'))\otimes V\cong SFH(X^\circ,\gamma^\circ|{G^{(m)}})\cong \widehat{HFK}(Y',K'|G^{(m)})\otimes V.$
\end{theorem}
\begin{proof}
    By Theorem 4.7, it suffices to prove that $SFH(M_0(K'),\gamma_0(K'))\cong QKFH'(Y',K',f).$

    Since all attaching curves in $\boldsymbol{\alpha}_{\text{left}}^{K'},\boldsymbol{\beta}_{\text{left}}^{K'}$ lie in $\operatorname{int}(\Sigma_{K'}),$ the inclusion $\Sigma_{K'}\hookrightarrow \widetilde{\Sigma}_{\max}^\circ$ identifies the intersection generators of the two diagrams.

    By our construction, $\widetilde{\Sigma}_{\max}^\circ =\Sigma_{K'}\cup _{\partial T} T,$ where $T\cong R_{\text{out}}$. Choose a one-holed torus $H\subset\operatorname{Int} T$. After cutting and capping, write
\[P=\overline{\widetilde{\Sigma}_{\max}^\circ\setminus H}\cup_{\partial H} D_p,\qquad
 \widehat{H}=H\cup_{\partial H}D_q\cong T^2;
\]
here $p,q$ denote the centers of the capping disks $D_p,D_q.$

All attaching curves lie on $P$, and $g(P)=g+d'-1$. 
Since a domain has no boundary inside the connected surface $T$,
$n_H(A)=n_T(A)$ for a domain $A$ of the Heegaard diagram.

\begin{lemma}
Choose a smooth family $\{j_s\}_{s\in[0,1]}$ of complex
structures on $P$, with endpoint values
$j_\beta:=j_0$ and $j_\alpha:=j_1$, such that $j_s=j_\beta$ near 0 and $j_s=j_\alpha$ near 1; also, on fixed neighborhoods of $p$ and the knot punctures,
$j_s$ is independent of $s$.

Fix a complex structure $j_{\widehat{H}}$ on $\widehat{H}$.
Plumb $(P,j_\beta)$ and $(P,j_\alpha)$ separately to $(\widehat{H},j_{\widehat{H}})$
at the cap centers $p,q$: remove small disks about these points
and join the resulting boundary circles by a complex cylinder
of length $\ell$. Denote the resulting manifold as $\widehat\Sigma_l$ and denote the resulting structures on
$\widehat\Sigma_\ell$ by
$j_{\beta,\ell}$ and $j_{\alpha,\ell}$.
We may arrange that the chocices of complex structures give us the following properties:

\begin{enumerate}
\item[(a)] \emph{\textbf{No positive-domain endpoint bubbles.}}
By varying the neck length $\ell>0$, we may choose a corresponding
family $\{J_\ell\}_{\ell>0}$ of admissible almost-complex structures, where each $J_\ell$ is defined on $\widehat\Sigma_\ell\times[0,1]\times\mathbb R$.

For some fixed $0<\epsilon_0<1/2$, require
\[
J_\ell=
\begin{cases}
j_{\beta,\ell}\oplus j_D, & 0\le s\le\epsilon_0,\\
j_{\alpha,\ell}\oplus j_D, & 1-\epsilon_0\le s\le1,
\end{cases}
\qquad\text{for every }\ell.
\]

Let $\mathcal A$ be the collection consisting of the nonnegative normalized domains
between the relevant generators with index in $\{-1,0,1,2\}$, and let $E_0$ denote the energy bound. There exists $l_0>0$ such that for any sequence $u_n$ of $J_{\ell_n}$-holomorphic
cylindrical curves satisfying
\[A(u_n)\in\mathcal A,\qquad
E(u_n)\le E_0,\qquad
\ell_n\ge\ell_0,
\]
no stable limit of this sequence contains a component
with positive surface domain whose strip projection is
a single point of $\{0,1\}\times\mathbb R$.

This applies both at fixed neck length and as
$\ell_n\to\infty$. The same exclusion holds for continuation
at a fixed neck length $\ell\ge\ell_0$, provided the
endpoint collars remain unchanged.

\item[(b)] \emph{\textbf{Uniform separation of visits to $\boldsymbol{p}$.}}
On the fixed surface $P$, use $J_P=j_s\oplus j_D$.
There exists $\epsilon_0<\delta<1/2$ such that for each cylindrical $J_P$-holomorphic curve $u=(u_P,u_D)$ with the usual attaching-curve boundary conditions, generator
ends, and domain in $\mathcal A_P$, we have
\[
V(s,t)\in D_p:=\{p\}+\operatorname{Sym}^{d'}(P)
\quad\Longrightarrow\quad
\delta\le s\le1-\delta;
\]
here, $\mathcal A_P$ consist of all nonnegative subdomains
of the $P$-side restrictions of $A\in\mathcal A$,
capped at $p$ with their neck multiplicities, $u_P$ and $u_D$ are the surface and strip projections,
respectively, and $V(s,t):=\sum_{w\in u_D^{-1}(s,t)} u_P(w)
\in\operatorname{Sym}^{d'+1}(P).$

Thus a sheet can pass through $p$ only at a strip height
uniformly separated from $0$ and $1$.
The same $\delta$ applies to every level of a Floer-breaking
limit; no embeddedness assumption on the associated incidence
curves is required. 
\end{enumerate}
\end{lemma}

\begin{proof}

In the discussion following the definition of $QKFH'$, we srtetch the neck along curves $c_\alpha^i$ to exclude boundary-degeneration. Since $H$ is disjoint from all attaching curves, we may assume each $c_i^\alpha$ are on the $P$-side. Then, similarly, we may srtetch the neck along curves $c_\alpha^i$ on $P$, and when the neck length $\ell'$ tends to infinity, $P$ converges to the nodal union $\widehat{H}_1^\alpha\cup \cdots\cup \widehat{H}_{d'}^\alpha\cup \widehat{R}'_\alpha$. Here, $\widehat{R}'_\alpha$ has genus $g-1>0$ and $\alpha_r$ remains essential on $\widehat{R}'_\alpha.$

The component of the limiting symmetric product containing
the attaching torus is therefore identified with
\[
X'_\alpha=\widehat H_1^\alpha\times\cdots\times
  \widehat H_{d'}^\alpha\times\widehat R'_\alpha.
\]
Now, let \[
L'_\alpha=\alpha_1\times\cdots\times\alpha_{d'}\times\alpha_r.\]
The same relative homotopy calculation as before gives
$\pi_2(X'_\alpha,L'_\alpha)=0$ and $\pi_2(X'_\alpha)=0$.
Hence every holomorphic disk in $X'_\alpha$ with boundary in $L'_\alpha$,
and every holomorphic sphere in $X'_\alpha$, has zero area and
is constant.

We now show that sufficiently long finite auxiliary necks
exclude nonconstant holomorphic disks satisfying the fixed domain restrictions and energy bound. Suppose otherwise; then choose a sequence of neck length $\ell'_n\to\infty$ and a corresponding sequence of
nonconstant disks
\[
v_n:(\mathbb D,\partial\mathbb D)
\longrightarrow
\bigl(\operatorname{Sym}^{d'+1}(P,j_{\alpha,\ell'_n}),
    \mathbb{T}_{\alpha^{K'}_{\text{left}}})
\]
within the prescribed bounds.
Choose a small neighborhood $U$ of the attaching torus,
identified throughout the family, which carries a symplectic form $\omega_U=d\theta$ that tames
the local complex structure, with $\theta$ vanishing on $\mathbb{T}_{\alpha^{K'}_{\text{left}}}$.
By Stokes' theorem, each $v_n$ must leave $U$.
Thus, for each $n$, we can choose a marked point $z_n$ with $v_n(z_n)\in\partial U$.

The domain bounds give uniform area bounds, and zero
multiplicity near the knot punctures prevents escape.
Marked compactness therefore gives a limiting map $v_\infty$ from
a connected collection of disks and spheres, with the
marked value still on $\partial U$.

Every boundary disk of this limit maps into $X'_\alpha$ with
boundary in $L'_\alpha$, so it is constant, say at $o\in L'_\alpha$.
Any disk or sphere attached to it must also pass through
$o$. Since $o$ contains no surface node, it belongs to
no other target component; the attached component
therefore lies in $X'_\alpha$ and is also constant at $o$.
Repeating along the connected source shows that the
entire limiting map has image $\{o\}$.
This contradicts the marked value lying on
$\partial U$, which is disjoint from $L'_\alpha$.

Thus all bounded nonconstant alpha-boundary disks are
excluded for sufficiently large finite $\ell'$. Fix a finite auxiliary neck length $\ell'_\alpha$ large
enough for the preceding exclusion, and let
$j_\alpha:=j_{\alpha,\ell'_\alpha}$ be the resulting
complex structure on $P$. Construct $j_\beta$ identically.

Keep the finite cylinder lengths $\ell'_\alpha,\ell'_\beta$ used in the above constructions fixed. Now, consider $\widehat\Sigma_\ell$.
As $\ell\to\infty$, the surface converges to
$P\cup_{p\sim q}\widehat H$. In the limiting symmetric product, the component containing
$\mathbb T_{\alpha_{\text{left}}^{K'}}$ is $\operatorname{Sym}^{d'+1}(P)$.
It meets the other components only along $D_p=\{p\}+\operatorname{Sym}^{d'}(P),$ where $D_p\cap\mathbb T_{\alpha_{\text{left}}^{K'}}=\varnothing.$

We have shown that every holomorphic disk in $(\operatorname{Sym}^{d'+1}(P),\mathbb T_{\alpha_{\text{left}}^{K'}})$
satisfying the prescribed domain and energy bounds is
constant. There are also no nonconstant holomorphic
spheres in $\operatorname{Sym}^{d'+1}(P)$. Indeed, a holomorphic sphere
\[
u:\mathbb{CP}^{1}\longrightarrow\operatorname{Sym}^{N}(P)
\]
has an incidence curve 
\[
\Gamma_u
=
\bigl\{(z,x)\in\mathbb{CP}^{1}\times\operatorname{int}P
:\ x\text{ occurs in the divisor }u(z)\bigr\}.
\]
whose normalization---the smooth compact
curve obtained by separating the local branches where divisor
points collide---maps holomorphically to $\operatorname{int}P$.
On each connected component, a nonconstant map would have open
image by the open mapping theorem. Its image is also compact,
hence closed; connectedness would then force it to equal
$\operatorname{int}P$, contradicting the noncompactness of
$\operatorname{int}P$. Thus every component map, and hence $u$,
is constant.

Suppose now that nonconstant alpha-boundary disks existed
for a sequence $\ell_n\to\infty$. As in preceding analysis, mark on each disk a point where it leaves
a fixed exact neighborhood of $\mathbb T_{\alpha_{{\text{left}}}}^{K'}$.
Compactness gives a stable limiting map in $\operatorname{Sym}^{d'+1}(P)$, with the marked value still outside that
neighborhood.

Choose a boundary-bearing disk component of the limit.
It lies in $\operatorname{Sym}^{d'+1}(P)$ and is constant
at some $o'\in\mathbb T_{\alpha_{{\text{left}}}}^{K'}$.
Every component attached to it must pass through $o'$.
Since $o'\notin D_p$, it cannot lie in another component
of the central fibre; it therefore remains in
$\operatorname{Sym}^{d'+1}(P)$ and is constant by the
preceding disk and sphere exclusions.
Repeating this argument over the connected limiting source
shows that every component is constant at $o'$.
This contradicts the marked value lying outside the chosen
neighborhood of $\mathbb T_{\alpha_{{\text{left}}}}^{K'}$.

Hence, alpha-boundary disk configurations are excluded for
all sufficiently large finite $\ell$, as well as at the
nodal limit. The beta case is identical.

 Write the original curves as $u_n=(u_{\Sigma,n},u_{D,n})$, with surface
and strip projections $u_{\Sigma,n},u_{D,n}$. On the fixed split alpha
collar, $u_{D,n}$ is a proper holomorphic branched covering of
degree $N=d'+1$. Hence
\[
V_n(z):=\sum_{w\in u_{D,n}^{-1}(z)}m_wu_{\Sigma,n}(w)
\]
defines a holomorphic map, where $m_w$ is the local degree of $u_{D,n}$ at $w$. Thus $V_n(z)$ is
the unordered degree-$N$ divisor of surface points lying
above $z$. Also,
\[
V_n(1,t)\in\mathbb T_{\alpha_{\text{left}}^{K'}} .
\]

Suppose a limiting component with positive surface domain
has constant strip projection at an alpha-boundary point.
After translating in $t$, take that point to be
$z_0=(1,0)$. Let $D^+$ be a fixed half-disk in the split alpha collar,
centered at $z_0$, and restrict $V_n$ to $D^+$.

Boundary compactness gives a stable limit of the maps
$V_n|_{D^+}$. Our previous endpoint-disk and sphere exclusions rule out
all nonconstant bubbles, while
$D_p\cap\mathbb T_{\alpha_{\text{left}}^{K'}}=\varnothing$ prevents an attached
component from entering another component of the limiting
target. Hence, after passing to a subsequence, $V_n|_{D^+}$
converges smoothly up to $z_0$. Let $V_\infty$ be the limit.

Choose a point $x_*$ in a region where the alleged
component has positive multiplicity, away from the
attaching curves and surface nodes. Then there are points
$w_n$ on the original curves such that
\[
u_{\Sigma,n}(w_n)\longrightarrow x_*,
\qquad
u_{D,n}(w_n)\longrightarrow z_0.
\]
Since $u_{\Sigma,n}(w_n)$ is contained in the divisor
$V_n(u_{D,n}(w_n))$, the limit divisor $V_\infty(z_0)$ contains
$x_*$. This is impossible because
\[
V_\infty(z_0)\in\mathbb T_{\alpha_{\text{left}}^{K'}},
\]
whereas $x_*$ lies on no alpha attaching curve.

The beta case is identical. The same local argument is
uniform for all sufficiently large neck lengths, including
the nodal limit, and for fixed-neck continuation preserving
the split endpoint collars. This proves (a). 

We now use the disk and sphere exclusions on $P$ established above to prove (b). Keep $P$ and $J_P=j_s\oplus j_D$ fixed;
no further neck stretching is involved.

For each cylindrical curve $u=(u_P,u_D)$ in the statement,
let $V$ be its associated divisor-valued map.
Near $s=0$, the map $V$ is holomorphic into
$\operatorname{Sym}^{d'+1}(P,j_\beta)$ and satisfies
$V(0,t)\in\mathbb T_\beta$.
The uniform area bound for domains in $\mathcal A_P$,
avoidance of the knot-puncture collars, and the preceding
disk and sphere exclusions give local compactness up to
this boundary without bubbling.

Suppose visits to $p$ are not uniformly separated from
$s=0$. Then there are curves $u_n=(u_{P,n},u_{D,n})$
as above, with associated maps $V_n$, and base points
$(s_n,t_n)\in[0,1]\times\mathbb R$ such that
\[
V_n(s_n,t_n)\in D_p,\qquad s_n\to0.
\]
Set
\[
\widetilde V_n(s,t):=V_n(s,t+t_n).
\]
On the beta collar, the translated maps $\widetilde V_n$
remain holomorphic into the fixed target
$\operatorname{Sym}^{d'+1}(P,j_\beta)$, satisfy
$\widetilde V_n(0,t)\in\mathbb T_{\beta_{\text{left}}^{K'}}$, and retain the
same area bound and puncture avoidance. The preceding compactness
argument therefore gives a subsequence converging smoothly
to $V_\infty$ on a fixed half-disk about $(0,0)$, up to the
beta boundary. Since $(s_n,0)\to(0,0)$,
\[
\widetilde V_n(s_n,0)=V_n(s_n,t_n)
\longrightarrow V_\infty(0,0).
\]
Every divisor on the left belongs to the closed set $D_p$,
so $V_\infty(0,0)\in D_p$.

On the other hand, $\widetilde V_n(0,t)\in\mathbb T_{\beta_{\text{left}}^{K'}}$,
and convergence up to the boundary gives
$V_\infty(0,0)\in\mathbb T_{\beta_{\text{left}}^{K'}}$.
Thus the same divisor belongs to both $D_p$ and
$\mathbb T_{\beta_{\text{left}}^{K'}}$, contradicting their disjointness.

The same argument at $s=1$ gives separation from the
alpha boundary. Taking the smaller separation constant,
we obtain $0<\delta<1/2$ such that
\[
V(s,t)\in D_p
\quad\Longrightarrow\quad
\delta\le s\le1-\delta.
\]
Every level of a Floer-breaking limit satisfies the same
domain and energy bounds, so the same estimate applies
to it. No embeddedness assumption on the associated
incidence curves is used.
\end{proof}
Let $\delta$ be the separation constant from
Lemma 4.9 (b). Choose $0<\epsilon_0\leq\epsilon<\delta$ within
the constant endpoint collars of $j_s$, and set
\[
I_\epsilon=[0,\epsilon]\cup[1-\epsilon,1].
\]
The subsequent plumbed family will use the common
split-collar width $\epsilon_0=\epsilon$.

\begin{lemma}
There exists a $t$-invariant almost-complex structure
$J_{\widehat{H}}$ on $\widehat{H}\times [0,1]_s\times \mathbb{R}_t$, satisfying $\boldsymbol{J1}$--$\boldsymbol{J4}$
and $\boldsymbol{J5'}$, which equals $j_{\widehat{H}}\oplus j_D$ near
$q$, on any fixed point neighborhoods required by
$\boldsymbol{J2}$, and whenever $s\in I_\epsilon$; for this almost complex structure, every nonconstant $J_{\widehat{H}}$-holomorphic map from a compact,
connected Riemann surface without boundary has image in
\[
\widehat{H}\times\{(s_0,t_0)\},
\qquad s_0\in I_\epsilon.
\]
This includes multiple covers.
\end{lemma}

\begin{proof}

Since $\widehat{H}$ is a torus, choose a smooth closed one-form
$\eta_0$ with $[\eta_0]\ne0$ in $H^1(\widehat{H};\mathbb R)$.
On small disjoint disks $U_i$ around $q$ and any
points on $\widehat{H}$ used in \textbf{J2}, the Poincar\'e lemma gives
functions $\theta_i$ with $\eta_0|_{U_i}=d\theta_i$.
Choose smooth cutoff functions $\chi_i$ supported in
$U_i$ and equal to $1$ on the smaller neighborhoods $U'_i$
where the split structure must be preserved. Then
\[
\eta:=\eta_0-\sum_i d(\chi_i \theta_i)
\]
vanishes on each $U'_i$, remains closed, and
satisfies $[\eta]=[\eta_0]\ne0$. Each $\chi_i \theta_i$ is
extended by zero outside $U_i$.

Choose a smooth function $\rho:[0,1]\to[0,1]$ that
vanishes on $I_\epsilon$ and is positive elsewhere.
Set $\lambda=\kappa\rho$, for a small constant $\kappa>0$,
and define
\[
\begin{aligned}
J_{\widehat{H}}\partial_s&=\partial_t,
& J_{\widehat{H}}\partial_t&=-\partial_s,\\
J_{\widehat{H}}X&=j_{\widehat{H}}X+\lambda(s)\eta(X)\partial_s
             +\lambda(s)\eta(j_{\widehat{H}}X)\partial_t,
&&X\in T\widehat{H}.
\end{aligned}
\]
The added terms vanish wherever $\eta=0$ or $\lambda=0$,
so all required split regions are preserved.

    These formulas give $J_{\widehat{H}}^2=-\mathrm{id}$.
For $W=X+a\partial_s+b\partial_t$,
\[
\Omega(W,J_{\widehat{H}}W)
=
\omega_{\widehat{H}}(X,j_{\widehat{H}}X)+a^2+b^2
+\lambda(s)\bigl(a\eta(j_{\widehat{H}}X)-b\eta(X)\bigr).
\]
The last term is uniformly small relative to the positive
quadratic part when $\kappa$ is small, so $\Omega$ tames
$J_{\widehat{H}}$. Moreover, with
\[
Y_X=X-\lambda(s)\eta(j_{\widehat{H}}X)\partial_s,
\]
we have
\[
J_{\widehat{H}}Y_X=Y_{j_{\widehat{H}}X},
\qquad
\Omega(Y_X,Y_Z)=\omega_{\widehat{H}}(X,Z).
\]
Thus $\xi=\{Y_X:X\in T\widehat{H}\}$ is an invariant symplectic
plane field as required by \textbf{J5'}.
The formula and the fixed split neighborhoods also give
\textbf{J1--J4}.

Let $u=(v,s,t):(C,j)\to \widehat{H}\times [0,1]\times \mathbb{R}$ be holomorphic,
where $C$ is compact, connected, and without boundary;
$v$ is its surface projection and $s,t$ its real-valued
strip coordinates. The equation $du\circ j=J_F\circ du$
gives
\[
dv\circ j=j_{\widehat{H}}\circ dv,
\qquad
ds\circ j=-dt+\lambda(s)v^*\eta.
\]
Define $\Lambda(r)=\int_0^r\lambda(\tau)\,d\tau$.
Since $d(v^*\eta)=0$,
\[
\begin{aligned}
ds\wedge(ds\circ j)
&=-ds\wedge dt+\lambda(s)\,ds\wedge v^*\eta\\
&=d\bigl(-s\,dt+\Lambda(s)v^*\eta\bigr).
\end{aligned}
\]
For any $j$-compatible metric on $C$, the left-hand side
is $-|ds|^2dA$. Stokes' theorem therefore gives
\[
\int_C|ds|^2dA=0,
\]
so $s\equiv s_0$.

The remaining equation is
\[
dt=\lambda(s_0)v^*\eta.
\]
If $v$ were constant, this would also make $t$ constant.
Thus, when $u$ is nonconstant, $v:C\to \widehat{H}$ is nonconstant.
It consequently has positive degree and induces an
injective pullback on real $H^1$, so $v^*[\eta]\ne0$.
Taking cohomology in the displayed equation gives
\[
0=[dt]=\lambda(s_0)v^*[\eta].
\]
Hence $\lambda(s_0)=0$, which means $s_0\in I_\epsilon$.
The same equation now gives $dt=0$, so $t\equiv t_0$.
Therefore the image of $u$ lies in the stated fibre.
\end{proof}

\begin{lemma}
Retain the choices of Lemmas 4.9 and 4.10. Let
$J_{\mathrm{aux},\ell}$ on
$\widehat\Sigma_\ell\times [0,1]\times \mathbb{R}=(P\#_\ell \widehat{H})\times [0,1]\times \mathbb{R}$
be obtained by plumbing $J_P=j_s\oplus j_D$ to $J_{\widehat{H}}$,
extending the mixed term by zero outside $H$.

There exists $\ell_*>0$ such that, for every
$\ell\ge\ell_*$, each embedded
$J_{\mathrm{aux},\ell}$-holomorphic cylindrical curve
$u:S\to\widehat\Sigma_\ell\times [0,1]\times \mathbb{R}$
of index one or two, with the prescribed boundary
conditions and generator ends, and satisfying
$A(u)\in\mathcal A$ and $E(u)\le E_0$, satisfies
\[
u(S)\cap(T\times [0,1]\times \mathbb{R})=\varnothing.
\]
\end{lemma}

\begin{proof}
The plumbed structures retain the admissibility conditions
and the common split endpoint collars.
Since $n_H(u)=n_T(u)$, positivity of intersections reduces
the desired avoidance of $T$ to proving $n_H(u)=0$.

Suppose this fails for arbitrarily large neck lengths.
Choose $\ell_n\to\infty$ and curves $u_n$ satisfying the
stated bounds with $n_H(u_n)>0$.
Since $\mathcal A$ is finite, we may pass to a subsequence of $u_n$ with a fixed surface domain $A$, where $n_H(A)=n_T(A)>0.$ Choose $x\in\operatorname{int}H$ away from the neck,
and for each curve in the subsequence, mark a point on the source whose surface projection
is $x$. Marked
cylindrical neck-stretching compactness gives a limiting
configuration, as in Section 12 of [5];
Lemma 4.9 (a) excludes positive-domain endpoint
components from this limit.

We now examine the $\widehat{H}$-side components of this limit.
To apply Lemma 4.10, we first show that they extend to
maps from closed Riemann surfaces.

All attaching curves and generator points lie on $P$, so the $\widehat H$-side components have no ordinary
boundary or generator ends. The compactness description
above leaves finitely many interior punctures; we show
that each is removable.

Near a puncture, choose a coordinate $\zeta=x+iy$ on
the source and write
\[
u=(v,s,t):\Delta^*\longrightarrow\widehat H\times [0,1]\times\mathbb{R},
\qquad \Delta^*=\Delta\setminus\{0\}.
\]
The surface projection $v$ is holomorphic and has finite
area by the domain bound. Since $\widehat H$ is compact,
$v$ extends holomorphically across $0$ by Proposition 5.16 of [12].

We next extend $h=s+it$. Since $v$ now extends smoothly,
\[
v^*\eta=a_1\,dx+a_2\,dy
\]
has smooth coefficients on the whole source disk.
The equation $ds\circ i=-dt+\lambda(s)v^*\eta$ gives
\[
\partial_{\bar\zeta}h
=\frac12\bigl(s_x-t_y+i(t_x+s_y)\bigr)
=\frac{\lambda(s)}2(-a_2+ia_1)=:\mathfrak{g}.
\]
The function $\mathfrak{g}$ is bounded near $0$, since
$0\le s\le1$. Assign it any value at $0$ and solve
\[
\partial_{\bar\zeta}\mathfrak{G}=\mathfrak{g},\qquad \mathfrak{G}\in W^{1,4}(\Delta),
\]
using the local $\bar\partial$ right inverse
in Proposition 2.70 of [12].
Sobolev embedding makes $\mathfrak{G}$ continuous and bounded.

Then $\mathfrak{f}:=h-\mathfrak{G}$ is holomorphic on $\Delta^*$ and
$\operatorname{Re}\mathfrak{f}=s-\operatorname{Re}\mathfrak{G}$ is bounded.
Hence $e^{\mathfrak{f}}$ and $e^{-\mathfrak{f}}$ extend holomorphically across $0$.
Their extensions still multiply to $1$, so the extension
of $e^{\mathfrak{f}}$ admits a local holomorphic logarithm.
This logarithm differs from $\mathfrak{f}$ by a constant in
$2\pi i\mathbb Z$, proving that $\mathfrak{f}$ extends.
Consequently, $h=\mathfrak{G}+\mathfrak{f}$ extends continuously and in
$W^{1,4}$; its displayed equation then gives smoothness
by elliptic bootstrapping.

Thus every interior puncture is removable, including
any alleged end with $t\to\pm\infty$.
At a surface-neck puncture, the extended surface value
is $q$. Filling the punctures and retaining ordinary
nodes as attachment points gives maps from closed
Riemann surfaces. The positive $H$-multiplicity retained
in the marked limit provides a nonconstant completed
component. Lemma 4.10 therefore gives its constant strip
value $(s_0,t_0)$, with $s_0\in I_\epsilon$.

We now show that the value $(s_0,t_0)$ also occurs
at a $P$-side evaluation in $D_p$.
The original source component giving rise to the marked
$\widehat H$-component has attaching boundary or generator
ends on $P$. Its connected limit therefore joins the
marked component to a $P$-side level carrying such
boundary or ends.

To show that the strip value $(s_0,t_0)$ is unchanged along this connection, start with the marked $\widehat H$-component
and stop at the first component with attaching boundary
or generator ends on $P$. After filling the source punctures, each of those components before
the last corresponds to a closed source. For source components mapping to
$\widehat H\times [0,1]\times \mathbb{R}$, Lemma 4.10 gives constant strip projection.
For source components mapping to the product of the stretched
cylinder with $[0,1]\times \mathbb{R}$, the strip projection extends holomorphically
to the closed source and is therefore constant.
Finally, a source component mapping to $P\times [0,1]\times \mathbb{R}$ has zero
symplectic area, since $H_2(P\times D;\mathbb R)=0$ and
the taming form is closed; tameness then forces the entire
map to be constant.

At each ordinary node or pair of matched surface-neck
punctures, the strip values on the two components agree.
Consequently, the value $(s_0,t_0)$ is preserved through
these attachments and is also the strip value at the
attachment to the $P$-side curve.

After filling its surface-neck punctures and including
trivial sheets, the $P$-side level gives a cylindrical
curve with the prescribed boundary conditions and domain
in $\mathcal A_P$. Its strip projection has degree
$N=d'+1$, inherited from the original attaching-boundary
sheets. Since the connection enters the $P$-side at $p$,
the associated divisor-valued map satisfies
\[
V(s_0,t_0)\in D_p.
\]

Lemma 4.9 (b) then gives
$s_0\in[\delta,1-\delta]$, whereas we have also shown that $s_0\in I_\epsilon$. This is impossible because
\[
I_\epsilon\cap[\delta,1-\delta]=\varnothing
\qquad(\epsilon<\delta).
\]
Consequently, $n_H(u)=n_T(u)=0$ for every curve in the
statement once $\ell$ is sufficiently large.
The positivity argument at the beginning then yields
\[
u(S)\cap(T\times [0,1]\times \mathbb{R})=\varnothing,
\]
as required.

\end{proof}
Fix a finite neck length $\ell_0$ for which Lemmas 4.9
and 4.11 apply, and write
$J_{\mathrm{aux}}=J_{\mathrm{aux},\ell_0}$.
Let $C$ be the common generator space and
$\partial_{\mathrm{aux}}$ the count of embedded
index-one auxiliary curves modulo translation.

By Lemma 4.11, every such curve avoids $T$, so its equation
agrees with the holomorphic-map equation defining the
sutured differential on $\Sigma_{K'}\times [0,1]\times \mathbb{R}$.
Conversely, every curve counted by the sutured differential
gives an auxiliary curve by extending its surface domain
by zero over $T$. This preserves the class index and
translation quotient. The linearizations of the
holomorphic-map equations also agree, since the
almost-complex structures coincide near each curve.
With the chosen regular data, the rigid moduli spaces
and their mod-two counts therefore agree. Hence
\[
(\mathcal{C},\partial_{\mathrm{aux}})
\cong
SFC
(\Sigma_{K'},\alpha_{\mathrm{left}}^{K'},
\beta_{\mathrm{left}}^{K'}),
\qquad
\partial_{\mathrm{aux}}^2=0.
\]

At the same neck length, choose regular split data
$J_{\mathrm{sp}}$ for $QKFH'(Y',K',f)$ with the same
endpoint and puncture models, and denote its count by
$\partial_{\mathrm{sp}}$.
The standing bounds and puncture avoidance give compactness
modulo breaking and bubbling. Lemma 4.9 (a) excludes
endpoint bubbling, while
$H_2(\widehat\Sigma_{\ell_0}\times [0,1]\times\mathbb{R};\mathbb R)=0$
and the closed taming form exclude nonconstant components
with closed source. The boundary count of the
one-dimensional index-two moduli spaces therefore gives
$\partial_{\mathrm{sp}}^2=0$.

Keep the neck length $\ell_0$ fixed. Let $J_0$ be the
split structure obtained from $J_{\mathrm{aux}}$ by
setting $\lambda=0$. Choose a smooth path $J^r$,
$r\in[0,1]$, from $J_{\mathrm{sp}}$ to
$J_{\mathrm{aux}}$: first interpolate through split
structures from $J_{\mathrm{sp}}$ to $J_0$, then
replace $\lambda$ by $m_\lambda\lambda$ as $m_\lambda$ increases from
$0$ to $1$. All structures agree on the fixed endpoint
collars and puncture neighborhoods.

Choose a smooth cutoff $\chi':\mathbb R\to[0,1]$ with
$\chi'(t)=0$ for $t\le-R$ and $\chi'(t)=1$ for $t\ge R$,
and set
\[
J_{\mathrm{cont}}(x,s,t)=J^{\chi'(t)}(x,s,t).
\]
Perturb this structure away from the fixed regions and
the two ends to achieve regularity, retaining the
notation $J_{\mathrm{cont}}$.

For generators $\mathbf x,\mathbf y$, let
$\mathcal M_{\mathrm{cont}}^0(\mathbf x,\mathbf y)$
denote the moduli space of embedded rigid
$J_{\mathrm{cont}}$-holomorphic curves from $\mathbf x$
to $\mathbf y$ in the allowed normalized classes of
combinatorial index zero. These curves satisfy the
usual attaching-curve boundary conditions. We quotient
by source reparametrization, but not by $t$-translation.
Define
\[
\Phi(\mathbf x)
=
\sum_{\mathbf y}
\#_2\mathcal M_{\mathrm{cont}}^0(\mathbf x,\mathbf y)
\,\mathbf y.
\]
The reverse interpolation defines $\Psi$.

The standing bounds and bubble exclusions give the
compactness required for these counts and their
one-parameter homotopies. The continuation and gluing
arguments in Section 9 of [5], with the corrected counting
conventions, show that $\Phi,\Psi$ are chain maps and
that their compositions are chain homotopic to the
identity. Hence
\[
(\mathcal{C},\partial_{\mathrm{sp}})
\simeq(\mathcal{C},\partial_{\mathrm{aux}}).
\]

Taking homology gives
\[
QKFH'(Y',K',f)
\cong H(C,\partial_{\mathrm{aux}})
\cong SFH(M_0(K'),\gamma_0(K')).
\]
Tensoring with $V$ and applying Theorem~4.7 and the
previously established knot Floer identification
completes the proof.\end{proof}

\section{Transfer of Floer Simplicity through Cut-and-Reglue}
\label{sec:transfer}

Combining the identifications from previous section with Corollary 20 in [7], we see that if $\operatorname{rank}\widehat{HF}(Y'|G^{(m)})=\operatorname{rank}\widehat{HFK}(Y',K'|G^{(m)})$ then $\operatorname{rank}SFH(M_0(K'),\gamma_0(K'))=2\operatorname{rank}SFH(M_0,\gamma_0),$ which is the rank condition needed to apply Binns' Theorem 1.1 [2].

When $K'$ is rationally null-homologous, we can just apply Ni's Proposition 6.8 [8] and the original Floer simple hypothesis to show that rank $\widehat{HF}(Y'|G)=$rank $\widehat{HF}(Y|G)=$rank $\widehat{HFK}(Y,K|G)=$rank $\widehat{HFK}(Y',K'|G).$ Also, as $K'$ is rationally null-homologous, rank $\widehat{HFK}(Y',K',\mathfrak{s})\geq$rank $\widehat{HF}(Y',\mathfrak{s})$ for any $\mathfrak{s}\in \mathrm{Spin}^c(Y')$. Thus, we see that rank $\widehat{HFK}(Y',K',\mathfrak{s})=$rank $\widehat{HF}(Y',\mathfrak{s})$ for any $\mathfrak{s}\in \mathcal{B}_{Y'}([G])$. Moreover, Ni's Lemma 7.4 [8] shows that for sufficiently large $m,$ $\mathcal{B}_{Y'}([G^{(m)}])\subset \mathcal{B}_{Y'}([G])$. All in all, $\operatorname{rank}\widehat{HF}(Y'|G^{(m)})=\operatorname{rank}\widehat{HFK}(Y',K'|G^{(m)})$.

The issue is that $K'$ need not be rationally null-homologous in $Y',$ so we cannot use the above argument directly.

The purpose of this section is to derive the rank equality without any
homological hypothesis on $K'$. We have generalized the definition of $\widehat{HFK}$ to knots that are not rationally null-homologous, we shall also generalize the surgery exact triangle used to prove Ni's Proposition 6.8 [8]. The ingredient is Alishahi--Eftekhary's surgery exact triangle for sutured
Floer homology [1], from which we derive a surgery triangle for meridian-sutured knot exteriors.

\subsection{A rank comparison by specialization}

Before applying the surgery triangle, we first establish a rank comparison that will be used in later argument. We consider three specializations of
the Alishahi--Eftekhary complex. Two recover the knot Floer and
Heegaard Floer groups, with the same puncturing factor; the third
provides a common lower bound for their dimensions. We show that
all three dimensions agree for the original pair \((Y,K)\).

Let \(J\subset W\) be a knot in a closed oriented three-manifold,
and work over \(\mathbb F=\mathbb F_2\). Equip $N_J=W\setminus\operatorname{int}(\nu J\sqcup B^3)$
with suture annuli \(\gamma_J\) around two oppositely oriented
meridians on \(\partial\nu J\) and an equator on \(\partial B^3\).

Cap the boundary circles of a sutured Heegaard diagram for
\((N_J,\gamma_J)\), placing basepoints \(w,z\) in the two
meridional capping disks and \(p\) in the remaining disk.
Choose the diagram to be admissible after forgetting \(z\).
In the refined sutured Floer chain complex associated with this
pointed diagram (See Section 5.3 of [1]), specialize the variables
corresponding to \(w,z,p\) by
\(u_w\mapsto0,\ u_z\mapsto u_z,\ u_p\mapsto0\).
Denote the resulting finite free \(\mathbb F[u_z]\)-complex by
\(C(W,J)\). Its differential counts Maslov-index-one holomorphic
disks with \(n_w=n_p=0\), weighted by \(u_z^{n_z}\).
Admissibility makes these sums finite, and Theorem 1.1 of [1]
implies that different admissible diagrams for \((N_J,\gamma_J)\)
give chain-homotopy equivalent \(\mathbb F[u_z]\)-complexes under
the above coefficient specialization.

Filling the \(z\)-suture gives \(W(2)\), the complement of two
balls in \(W\). For each generator \(\mathbf x\), let
\(\mathfrak{s}(\mathbf x)\in\operatorname{Spin}^c(N_J,\gamma_J)\)
be its relative structure, and denote by
\([\mathfrak{s}(\mathbf x)]\in\operatorname{Spin}^c(W)\) its extension
across this filling and the two capping balls.

For a disk \(\phi\) counted by the differential,
Corollary~2.9 of [1] gives
\begin{equation}\label{eq:weighted-label-difference}
\mathfrak{s}(\mathbf x)-\mathfrak{s}(\mathbf y)
=n_z(\phi)\operatorname{PD}[\mu_J],
\end{equation}
where \(\mu_J\) is oriented as the \(z\)-suture.
Since the meridian bounds a disk after filling,
\([\mathfrak{s}(\mathbf x)]=[\mathfrak{s}(\mathbf y)]\).
Thus, for \(s\in\operatorname{Spin}^c(W)\), the generators
with \([\mathfrak{s}(\mathbf x)]=\mathfrak{s}\) span a subcomplex
\(C_\mathfrak{s}=C(W,J,\mathfrak{s})\) over \(\mathbb F[u_z]\).

For \(a=0,1\), set \(C_\mathfrak{s}(a)=C_\mathfrak{s}/(u_z-a)C_\mathfrak{s}\).
At \(u_z=0\), all three basepoint multiplicities must vanish,
giving the sutured Floer complex of \((N_J,\gamma_J)\).
At \(u_z=1\), the \(z\)-weight disappears, giving the complex
of \(W(2)\). The knot-exterior identification and puncturing
formula [3] therefore yield
\begin{equation}\label{eq:coefficient-specializations}
\begin{aligned}
H_*(C_\mathfrak{s}(0))
&\cong\widehat{HFK}(W,J,\mathfrak{s})\otimes V,\\
H_*(C_\mathfrak{s}(1))
&\cong\widehat{HF}(W,\mathfrak{s})\otimes V,
\end{aligned}
\end{equation}
where \(V=\mathbb F^2\), with the same ambient label in both
identifications. Finally, define
\[
r(W,J;\mathfrak{s})
=
\dim_{\mathbb F(u_z)}
H_*\!\left(C_\mathfrak{s}\otimes_{\mathbb F[u_z]}\mathbb F(u_z)\right).
\]

\begin{lemma}\label{lem:coefficient-rank-comparison}
For every knot \(J\subset W\) and
\(\mathfrak{s}\in\operatorname{Spin}^c(W)\),
\begin{equation}\label{eq:generic-rank-bounds}
\begin{aligned}
r(W,J;\mathfrak{s})&\le 2\operatorname{rank}\widehat{HFK}(W,J,\mathfrak{s}),\\
r(W,J;\mathfrak{s})&\le 2\operatorname{rank}\widehat{HF}(W,s).
\end{aligned}
\end{equation}
If \(J\) is rationally null-homologous, then
\begin{equation}\label{eq:rational-generic-rank}
r(W,J;\mathfrak{s})=2\operatorname{rank}\widehat{HF}(W,\mathfrak{s}).
\end{equation}
\end{lemma}

\begin{proof}
If \(C_\mathfrak{s}=0\), the claims are immediate. Otherwise, let \(\mathbf{D}(u_z)\)
be its differential matrix and \(|C_\mathfrak{s}|\) its number of generators.

We first prove the inequalities. If \(\mathbf{D}(a)\) has rank \(n\), some \(n\times n\) submatrix of it has nonzero determinant. The determinat of the corresponding submatrix of $\mathbf{D}(u_z)$ is therefore a
nonzero polynomial, so \(\mathbf{D}(u_z)\) has rank at least \(n\) over
\(\mathbb F(u_z)\). Since
\(\dim H=\dim\ker \mathbf{D}-\dim\operatorname{im}\mathbf{D}
=|C_\mathfrak{s}|-2\operatorname{rank}\mathbf{D}\), we obtain, for \(a=0,1\),
\[
\begin{aligned}
r(W,J;\mathfrak{s})
&=|C_\mathfrak{s}|-2\operatorname{rank}_{\mathbb F(u_z)}\mathbf{D}(u_z)\\
&\le |C_\mathfrak{s}|-2\operatorname{rank}_{\mathbb F}\mathbf{D}(a)
=\dim_{\mathbb F}H_*(C_\mathfrak{s}(a)).
\end{aligned}
\]
The identifications in \eqref{eq:coefficient-specializations}
give the two asserted inequalities.

Now suppose that \(J\) is rationally null-homologous.
The class \(\kappa=\operatorname{PD}[\mu_J]\) has infinite order,
since a rational Seifert surface intersects the meridian
nontrivially. All generators of \(C_\mathfrak{s}\) have the same ambient
label, so their relative structures differ by multiples of
\(\kappa\) (See Section 2.2 of [1]). Fix a reference generator
\(\mathbf x_0\) and define \(h(\mathbf x)\in\mathbb Z\)
for every generator by
\[
\mathfrak{s}(\mathbf x)-\mathfrak{s}(\mathbf x_0)=h(\mathbf x)\kappa.
\]
These integers depend only on the generators. For every counted
disk \(\phi\) from \(\mathbf x\) to \(\mathbf y\),
\eqref{eq:weighted-label-difference} gives
\[
n_z(\phi)\kappa
=\bigl(h(\mathbf x)-h(\mathbf y)\bigr)\kappa.
\]
Since \(\kappa\) has infinite order,
\(n_z(\phi)=h(\mathbf x)-h(\mathbf y)\).
Thus all disk weights are determined by the same function \(h\).

Over \(\mathbb F(u_z)\), choose the basis
\(\mathbf e_{\mathbf x}=t^{-h(\mathbf x)}\mathbf x\).
For every counted disk, its contribution in this basis satisfies
\[
u_z^{-h(\mathbf x)}u_z^{n_z(\phi)}\mathbf y
=
u_z^{-h(\mathbf x)+n_z(\phi)+h(\mathbf y)}
\mathbf e_{\mathbf y}
=
\mathbf e_{\mathbf y}.
\]
Hence this single change of basis removes all \(u_z\)-weights
simultaneously, leaving exactly the differential matrix \(\mathbf{D}(1)\).
Field extension preserves matrix rank, so the total homology
dimensions agree:
\[
r(W,J;\mathfrak{s})
=\dim_{\mathbb F}H_*(C_\mathfrak{s}(1))
=2\operatorname{rank}\widehat{HF}(W,\mathfrak{s}).
\]
\end{proof}

\begin{corollary}\label{cor:initial-three-ranks}
For the original pair \((Y,K)\) and every
\(\mathfrak{s}\in\operatorname{Spin}^c(Y)\),
\begin{equation}\label{eq:initial-three-ranks}
r(Y,K;\mathfrak{s})
=
2\operatorname{rank}\widehat{HFK}(Y,K,\mathfrak{s})
=
2\operatorname{rank}\widehat{HF}(Y,\mathfrak{s}).
\end{equation}
\end{corollary}

\begin{proof}
Since \(K\) is rationally null-homologous,
Lemma~\ref{lem:coefficient-rank-comparison} gives
\[
2\operatorname{rank}\widehat{HFK}(Y,K,\mathfrak{s})
\ge r(Y,K;\mathfrak{s})
=2\operatorname{rank}\widehat{HF}(Y,\mathfrak{s}).
\]
Floer simplicity gives equality of the total knot Floer and
Heegaard Floer ranks. The nonnegative differences therefore
vanish in every ambient summand.
\end{proof}

We next use the surgery triangle to preserve the totals of
these three ranks in the \(G\)-bottommost degree. The inequalities
of Lemma~\ref{lem:coefficient-rank-comparison} will then allow us
to recover equality in each ambient summand of the reglued pair.

\subsection{Cut-and-reglue invariance}

For the auxiliary rank, write
\[
r(W,J\mid G)
=
\sum_{\mathfrak{s}\in B_W([G])}r(W,J;\mathfrak{s}).
\]
Only finitely many terms are nonzero.

\begin{lemma}
\label{lem:regluing-three-ranks}
For the pairs \((Y,K)\) and \((Y',K')\) related by cutting
and regluing along \(G\),
\begin{equation}\label{eq:regluing-three-ranks}
\begin{aligned}
\operatorname{rank}\widehat{HFK}(Y,K\mid G)
&=\operatorname{rank}\widehat{HFK}(Y',K'\mid G),\\
\operatorname{rank}\widehat{HF}(Y\mid G)
&=\operatorname{rank}\widehat{HF}(Y'\mid G),\\
r(Y,K\mid G)&=r(Y',K'\mid G).
\end{aligned}
\end{equation}
\end{lemma}

\begin{proof}
Since the gluing maps defining \(Y\) and \(Y'\) differ,
up to isotopy, by a composition of Dehn twists about
nonseparating curves in \(G\), it suffices to prove the claim
when \((Y',K')\) is obtained from \((Y,K)\) by cutting along
\(G\) and regluing by a single such twist.

Let \(c\subset G\) be the curve defining the chosen Dehn twist.
Push \(c\) slightly off \(G\), keeping it disjoint from \(K\),
and retain the notation \(c\). Let \(Y_r\) be obtained from \(Y\) by
Dehn surgery on \(c\) with coefficient \(r\), measured
relative to the framing induced by \(G\), and let
\(K_r\subset Y_r\) be the image of \(K\) in \(Y_r\).
The coefficient \(\infty\) denotes meridional filling,
so \((Y_\infty,K_\infty)=(Y,K)\).
The surgery description of the chosen Dehn twist identifies
\((Y_\varepsilon,K_\varepsilon)\) with \((Y',K')\) for
an appropriate sign \(\varepsilon=\pm1\).
For the surgery triangle, use the coefficients
\(\infty,0,\varepsilon\). Denote them by \(r_0,r_1,r_2\)
in the cyclic order of Section 8 in [1], choosing \(r_1=0\).

Set
\(
Z=Y\setminus\operatorname{int}(\nu K\cup\nu c).
\)
On \(\partial\nu K\), keep two oppositely oriented meridional
sutures, with coefficient variables \(u_w,u_z\).
For term \(i\), equip \(\partial\nu c\) with two annular
sutures \(\gamma_{i,1},\gamma_{i,2}\) whose cores are
oppositely oriented parallel curves representing the
filling slope determined by \(r_i\), using the framing
induced by \(G\). Denote their coefficient variables by
\(\zeta_{i,1},\zeta_{i,2}\).

Choose an admissible surgery Heegaard quadruple as in
Section 8.1 of [1]. It has marked points \(z_0,z_1,z_2\) lying beside the three sides of the small triangle formed by
the surgery curves, and an auxiliary point \(p\) inside the triangle, as in Figure 3 of [1]. Denote the variables representing these points in the
common coefficient algebra by
\(\lambda_0,\lambda_1,\lambda_2,\lambda_p\).
The coefficient algebra in Section 8.1 of [1] imposes the relation
\[
\lambda_p=\prod_{i=0}^{2}\lambda_i^{m_i-1},
\]
where \(m_i\) is the algebraic intersection number of the
surgery curves for terms \(i+1\) and \(i+2\), with indices
modulo three.
Here \(m_0=m_1=m_2=1\), so
\(\lambda_p=1\).
The coefficient embeddings therefore send
\[
\zeta_{i,1}\longmapsto\lambda_i,\qquad
\zeta_{i,2}\longmapsto\lambda_{i+1}\lambda_{i+2},
\]
with indices modulo three, and leave \(u_w,u_z\) unchanged.

By Section 8.4 of [1], the specialization
\(\lambda_0=\lambda_1=1\) corresponds to filling
\(\gamma_{0,1}\), \(\gamma_{1,1}\), and \(\gamma_{2,2}\),
respectively, leaving the variable \(\lambda_2\) on each
surviving surgery suture. Each filling attaches
\(D^2\times I\) along the chosen suture annulus, replacing
the surgery torus by a sphere with one suture.
Capping this sphere completes the corresponding Dehn filling;
hence, before capping, the manifold is
\(
Y_{r_i}\setminus\operatorname{int}(\nu K_{r_i}\cup B^3).
\)

Put \(R=\mathbb F[u_z]\). The coefficient algebra of each
filled sutured manifold is
\(\mathbb F[\lambda_2,u_w,u_z]\), since its positive and
negative boundary regions are each an annulus and a disk
(See Section 3.4 of [1]). We therefore apply the homomorphism
\[
\mathbb F[\lambda_2,u_w,u_z]\longrightarrow R,
\qquad
\lambda_2\longmapsto0,\quad
u_w\longmapsto0,\quad
u_z\longmapsto t.
\]
This gives the \(R\)-complexes
\[
C_i=C(Y_{r_i},K_{r_i})
=\bigoplus_{\mathfrak{s}\in\operatorname{Spin}^c(Y_{r_i})}
C(Y_{r_i},K_{r_i},\mathfrak{s}),
\qquad i=0,1,2.
\]

Apply the same coefficient specialization to the chain maps defined by counting holomorphic triangles in Section~8.1 of [1]. In their construction, generators are
grouped according to the \(\operatorname{Spin}^c\) structure
to which their relative structures extend after all sutures
are filled. Taking the direct sum of the specialized maps
between these summands gives chain maps \(C_i\to C_{i+1}\),
with indices modulo three. By Theorem~8.6 and Section~8.4 of [1], these
maps induce the following exact triangle after tensoring
with any \(B\in\{R/(u_z),R/(u_z-1),\mathbb F(u_z)\}\):
\begin{equation}\label{eq:coefficient-surgery-triangle}
\begin{tikzcd}[column sep=small, row sep=large]
H_*(C_0\otimes_R B) \arrow[rr]
&& H_*(C_1\otimes_R B) \arrow[dl]\\
& H_*(C_2\otimes_R B) \arrow[ul]&
\end{tikzcd}
\end{equation}

By Section~5.1, the terms satisfy
\[
\begin{aligned}
H_*(C_i\otimes_R R/(u_z))
&\cong \widehat{HFK}(Y_{r_i},K_{r_i})\otimes V,\\
H_*(C_i\otimes_R R/(u_z-1))
&\cong \widehat{HF}(Y_{r_i})\otimes V,
\end{aligned}
\]
and
\[
\dim_{\mathbb F(u_z)}H_*(C_i\otimes_R\mathbb F(u_z))
=
\sum_{\mathfrak{s}\in\operatorname{Spin}^c(Y_{r_i})}
r(Y_{r_i},K_{r_i};\mathfrak{s}).
\]

We next show that the triangles preserve the Chern evaluation
defining the \(G\)-bottommost summands.

Choose a parallel copy \(G_\delta\subset\operatorname{int}Z\)
of \(G\). For \(i=0,2\), let \(G_i\subset Y_{r_i}\) be its
image under the filling inclusion \(Z\hookrightarrow Y_{r_i}\).
In \(Y_{r_1}=Y_0\), compress its image along the curve
parallel to the original twisting curve, using the disk
provided by the zero filling, and denote the resulting
surface by \(G_1\). Since the twisting curve is nonseparating,
\[
g(G_0)=g(G_2)=g(G),\qquad g(G_1)=g(G)-1.
\]

Compression preserves the homology class, so \([G_\delta]\)
maps to \([G_i]\) in every filling. The product
\(G_\delta\times[0,1]\) identifies these classes in each
surgery cobordism, so the holomorphic triangle maps preserve
the integer pairing \(\langle c_1(\mathfrak s),[G_i]\rangle\). Moreover, multiplication by coefficient monomials shifts relative
\(\operatorname{Spin}^c\) labels by sums of Poincar\'e duals
of curves on \(\partial Z\). These shifts pair trivially
with \([G_\delta]\), so the
chain maps and chain homotopies in A--E's proof of exactness preserve the Chern evaluation after
specialization (See Theorem~8.6 of [1]). The triangles therefore
restrict to the summands with
\(\langle c_1(\mathfrak s),[G_i]\rangle=2-2g(G)\).

In the restricted triangles above, every zero-surgery summand satisfies
\[
|\langle c_1(\mathfrak s),[G_1]\rangle|
=2g(G)-2>2g(G_1)-2.
\]
Since \(G_1\) is disjoint from \(K_0\), the adjunction inequalities ([9], [11]) and Lemma 5.1 give
\[
\widehat{HFK}(Y_0,K_0,\mathfrak s)
=\widehat{HF}(Y_0,\mathfrak s)=0,
\qquad r(Y_0,K_0;\mathfrak s)=0.
\]
Thus the zero-surgery homology vanishes for every coefficient choice, and exactness gives an isomorphism between the restricted homologies of the \(\infty\)-filling $(Y,K)$ and the \(\varepsilon\)-filling $(Y',K')$.

Since \(\chi(G_0)=\chi(G_2)=2-2g(G)\), these restrictions select exactly the \(G\)-bottommost summands. Applying (4), which preserves ambient \(\operatorname{Spin}^c\) labels, gives
\[
\begin{aligned}
\widehat{HFK}(Y,K\mid G)\otimes V
&\cong \widehat{HFK}(Y',K'\mid G)\otimes V,\\
\widehat{HF}(Y\mid G)\otimes V
&\cong \widehat{HF}(Y'\mid G)\otimes V.
\end{aligned}
\]
Over \(\mathbb F(u_z)\), we similarly obtain \(r(Y,K\mid G)=r(Y',K'\mid G)\).

We have thus proven the lemma for regluing by a single Dehn twist, and the general case follows by iteration.
\end{proof}

\begin{proposition}
\label{prop:final-bottommost-equality}
For every \(\mathfrak{s}'\in B_{Y'}([G])\),
\begin{equation}\label{eq:final-bottommost-equality}
\operatorname{rank}\widehat{HFK}(Y',K',\mathfrak s')
=
\operatorname{rank}\widehat{HF}(Y',\mathfrak s').
\end{equation}
\end{proposition}

\begin{proof}
Section~5.1 gives equality of the three ranks for the original
pair. Summing over \(B_Y([G])\) and applying
Lemma~\ref{lem:regluing-three-ranks}, we obtain
\begin{equation}\label{eq:final-three-totals}
r(Y',K'\mid G)
=
2\operatorname{rank}\widehat{HFK}(Y',K'\mid G)
=
2\operatorname{rank}\widehat{HF}(Y'\mid G).
\end{equation}

For each \(\mathfrak s'\in B_{Y'}([G])\), Lemma~5.1 gives
\[
\begin{aligned}
2\operatorname{rank}\widehat{HFK}(Y',K',\mathfrak s')
-r(Y',K';\mathfrak s')&\ge0,\\
2\operatorname{rank}\widehat{HF}(Y',\mathfrak s')
-r(Y',K';\mathfrak s')&\ge0.
\end{aligned}
\]
By \eqref{eq:final-three-totals}, the sum of each expression
over \(B_{Y'}([G])\) is zero. Every term therefore vanishes,
giving
\[
2\operatorname{rank}\widehat{HFK}(Y',K',\mathfrak s')
=
r(Y',K';\mathfrak s')
=
2\operatorname{rank}\widehat{HF}(Y',\mathfrak s').
\]
Cancelling \(2\) proves the claim.
\end{proof}

\subsection{The $G^{(m)}$-bottommost summands}
\label{subsec:Gm-bottommost}

We now apply Proposition~5.4 to the summands appearing in the closure
comparison of Section~4. Fix the reglued pair $(Y',K')$. The oriented cut-and-paste construction in Section 3 gives
\begin{equation}
\label{eq:Gm-affine-data}
[G^{(m)}]=[\overline S]+m[G],
\qquad
\chi(G^{(m)})=\chi(\overline S)+m\chi(G).
\end{equation}
Here $(Y',K')$ and $\overline S$ are independent of $m$, and
$\chi(\overline S)=\chi(S)$. By Section~3, $G^{(m)}$ is connected;
since $\chi(G)<0$, it has genus greater than one for all sufficiently
large $m$.

Let
\[
\mathcal S=
\bigl\{\mathfrak s'\in\operatorname{Spin}^c(Y'):
\widehat{HFK}(Y',K',\mathfrak s')\ne0
\ \text{or}\
\widehat{HF}(Y',\mathfrak s')\ne0\bigr\}.
\]
Only finitely many ambient labels occur among the generators of
$C(Y',K')$, so the identifications~(4) imply that $\mathcal S$ is finite.
The following lemma uses the finite-set argument of Lemma 7.4 in [8].

\begin{lemma}
\label{lem:Gm-supported-selection}
There is an integer $m_0>0$ such that, for every $m\ge m_0$,
\[
\mathcal S\cap\mathcal B_{Y'}([G^{(m)}])
\subseteq\mathcal B_{Y'}([G]).
\]
\end{lemma}

\begin{proof}
Choose $m_0$, also satisfying the largeness requirements of Section~3,
such that
\[
\bigl|\langle c_1(\mathfrak s'),[\overline S]\rangle
-\chi(\overline S)\bigr|<m_0
\qquad\text{for every }\mathfrak s'\in\mathcal S.
\]
For $m\ge m_0$ and
$\mathfrak s'\in\mathcal S\cap\mathcal B_{Y'}([G^{(m)}])$,
the defining equality
$\langle c_1(\mathfrak s'),[G^{(m)}]\rangle=\chi(G^{(m)})$
and \eqref{eq:Gm-affine-data} give
\[
m\bigl(\langle c_1(\mathfrak s'),[G]\rangle-\chi(G)\bigr)
=\chi(\overline S)-\langle c_1(\mathfrak s'),[\overline S]\rangle.
\]
The right-hand side has absolute value less than $m$. Thus the integer
$\langle c_1(\mathfrak s'),[G]\rangle-\chi(G)$ has absolute value
less than one and must vanish. Hence
$\mathfrak s'\in\mathcal B_{Y'}([G])$.
\end{proof}

\begin{corollary}
\label{cor:Gm-bottommost-rank-equality}
For all sufficiently large $m$,
\begin{equation}
\label{eq:Gm-bottommost-rank-equality}
\operatorname{rank}\widehat{HFK}(Y',K'\mid G^{(m)})
=\operatorname{rank}\widehat{HF}(Y'\mid G^{(m)}).
\end{equation}
\end{corollary}

\begin{proof}
Fix $m\ge m_0$ and
$\mathfrak s'\in\mathcal B_{Y'}([G^{(m)}])$.
If $\mathfrak s'\notin\mathcal S$, both
$\widehat{HFK}(Y',K',\mathfrak s')$ and
$\widehat{HF}(Y',\mathfrak s')$ vanish. Otherwise,
Lemma~\ref{lem:Gm-supported-selection} gives
$\mathfrak s'\in\mathcal B_{Y'}([G])$, and Proposition~5.4 gives
\[
\operatorname{rank}\widehat{HFK}(Y',K',\mathfrak s')
=\operatorname{rank}\widehat{HF}(Y',\mathfrak s').
\]
Summing over $\mathcal B_{Y'}([G^{(m)}])$ proves the assertion;
only finitely many summands are nonzero.
\end{proof}

\section{Proof of Main Theorem}

We now assemble the preceding sections.  Fix $m$ large enough that Corollary~5.6 apply, so that
\[
\operatorname{rank}\widehat{HFK}(Y',K'\,|\,G^{(m)})
=
\operatorname{rank}\widehat{HF}(Y'\,|\,G^{(m)}).
\]
By the closure comparison of Section 4, this equality is equivalent to
$$\operatorname{rank }SFH(M_0(K'),\gamma_0(K'))=2\operatorname{ rank }SFH(M_0,\gamma_0).$$
Moreover, $(M_{0}(K'),\gamma_{0}(K'))$ is a taut balanced sutured
manifold by discussion of Section 3, which means the group $SFH(M_{0}(K'),\gamma_{0}(K'))$ is
nonzero by Juh\'asz's nonvanishing theorem [4].
Hence the common
value in the rank equality above is nonzero.

It remains to verify Binns's homological hypothesis.
Let
\(i:\partial M_0\hookrightarrow M_0\) be the boundary inclusion. By the long exact sequence of the pair \((M_0,\partial M_0)\), the condition
\(
[K']\in\operatorname{im}i_*
\)
is equivalent to the image of \([K']\) vanishing in
\(H_1(M_0,\partial M_0;\mathbb Q)\). By Poincaré–Lefschetz duality, it therefore suffices to show that
\(
\mathcal{F}\cdot K'=0
\)
for every closed oriented surface
\(\mathcal{F}\subset\operatorname{int}M_0\).

Regard \(\mathcal{F}\) as a surface in \(Y'\), and choose one of the parallel copies \(G_j\) used to construct \(G^{(m)}\). Outside a small neighborhood of \(S\cap G_j\), the surface \(G_j\) is contained in \(G^{(m)}\). Since \(\mathcal{F}\cap G^{(m)}=\varnothing\), after making \(\mathcal{F}\) transverse to \(G_j\), the cycle \(\mathcal{F}\cap G_j\) is supported in this neighborhood and hence is homologous to a rational multiple of the common class of the components of \(S\cap G_j\). Cutting \(Y'\) open along \(G_j\) identifies the resulting manifold with \(M\), and cuts \(\mathcal{F}\) to a surface \(\mathcal{F}_M\subset M\) with
\(
[\partial\mathcal{F}_M]
=
q\bigl([c_+]-[c_-]\bigr)
\)
for some \(q\in\mathbb Q\).

By Section 3, \(
[\partial S]
=
k\bigl([c_+]-[c_-]\bigr).
\)
Thus the boundary of \(\mathcal{F}_M-\frac qk S\) is null-homologous in \(G_+\sqcup G_-=\partial M\), and it can be capped off in \(\partial M\) to give a closed rational two-cycle \(Z\subset M\). Cutting along \(G_j\) is disjoint from the knot, and both \(S\) and the added boundary chain are disjoint from \(K\). Hence
\(
\mathcal{F}\cdot K'
=
\mathcal{F}_M\cdot K
=
Z\cdot K
=
0,
\)
where the final equality follows from
\([K]=0\in H_1(Y;\mathbb Q)\). This proves Binns’s homological hypothesis.

Together with the nonzero rank equality established above,
this verifies all the hypotheses of Binns's Theorem 1.1 [2].
Consequently, $M_0$ has a closed connected summand $M_{\mathrm{cl}}$ containing $K'$, in which \(K'\) is Floer simple. Since \(M_0\) is irreducible, the corresponding connected-sum sphere bounds a ball. The side meeting \(\partial M_0\) has additional boundary and cannot be this ball, so the punctured \(M_{\mathrm{cl}}\)-side is a ball and
\(
M_{\mathrm{cl}}\cong S^3.
\)
Unknot detection [9] therefore shows that \(K'\) bounds an embedded disk \(D'\subset M_0\).

Choose \(D'\) transverse to the regluing surface \(G\), with
\(|D'\cap G|\) minimal. If \(D'\cap G\neq\varnothing\), choose a component \(\gamma\) innermost on \(D'\), and let \(\Delta\subset D'\) be the disk it bounds. Since \(\operatorname{int}\Delta\cap G=\varnothing\), incompressibility of \(G_\pm\subset M\) implies that \(\gamma\) also bounds a disk in \(G\). Choose such a disk \(E\subset G\) with \(\operatorname{int}E\cap D'=\varnothing\). Replacing the subdisk of \(D'\) bounded by \(\gamma\) with a slight pushoff of \(E\) reduces \(|D'\cap G|\), a contradiction.

Thus $D'$ is disjoint from $G$ and lies in
$Y'\setminus\!\setminus G\cong M$.
Under the original gluing $M\to Y$, this disk remains
embedded and its boundary becomes $K$.
Hence $K$ is unknotted. This completes the proof of Theorem 1.1.
 \qed

\noindent\section{Acknowledgment}

I would like to thank my advisor Professor Yi Ni for his guidance, and I would like to thank Alan Du for discussing the problem with me.

\end{document}